\documentclass[12pt]{article}

\usepackage{lscape}
\usepackage[english]{babel}
\usepackage{amssymb,amsmath}
\usepackage{graphicx,epsfig,color}
\usepackage{amsfonts,latexsym}
\usepackage{longtable}

\newtheorem{theorem}{Theorem}

\newtheorem{proposition}{Proposition}
\newtheorem{remark}{Remark}
\newtheorem{corollary}{Corollary}

\begin{document}
	%
	%

	\noindent GOODNESS-OF-FIT TESTS FOR THE BIVARIATE POISSON DISTRIBUTION
	\vskip 3mm

	\vskip 5mm
	\noindent F. Novoa-Mu\~noz
	
	\noindent Departamento de Estad\'istica
	
	\noindent Universidad del B\'io-B\'io
	
	
	\noindent fnovoa@ubiobio.cl
	
	\vskip 3mm

	\noindent ABSTRACT
	
	The bivariate Poisson distribution is commonly used to model bivariate count data. In this paper we study a goodness-of-fit test for this distribution. We also provide a review of the existing tests for the bivariate Poisson distribution, and its multivariate extension. The proposed test is consistent against any fixed alternative. It is also able to detect local alternatives converging to the null at the rate $n^{-\frac{1}{2}}$. The bootstrap can be employed to consistently estimate the null distribution of the test statistic. Through a simulation study we investigated the goodness of the bootstrap approximation and the power for finite sample sizes.
	
	\vskip 4mm

	\section{Introduction}\label{Introduction}
	The univariate Poisson distribution (UPD) has helped to model many real life situations. For a survey of statistical issues, problems and applications associated with the UPD the reader is referred to the text of Haight (1967) and Johnson and Kotz (1969). For the other hand, the bivariate Poisson distribution (BPD) is appropriate for modelling paired count data exhibiting positive correlation.
	
	Several definitions for the BPD have been given (see, e.g. Kocherlakota and Kocherlakota, 1992). In this paper we will work with the following one, because it has received the most attention in the statistical literature (see, e.g. Holgate, 1964; Johnson, Kotz and Balakrishnan, 1997). Let
	\[X_1=Y_1+Y_{3}\, \quad \textnormal{ and } \quad X_2=Y_2+Y_{3},\]
	where $Y_1,Y_2$ and $Y_{3}$ are independent Poisson random variables with means $\theta'_1=\theta_1-\theta_{3}>0$, $\theta'_2=\theta_2-\theta_{3}>0$ and $\theta_{3}>0$, respectively. The joint distribution of the vector $(X_1,  X_2)$ is called BPD with parameter $\theta=(\theta_1,\theta_2,\theta_3)$, $(X_1, X_2)\sim BP(\theta)$ for short.
	
	In the statistical literature on goodness-of-fit (gof) tests for the BPD, which is rather sparse in comparison with the univariate case, we found the following: the tests given by Crockett (1979), Loukas and Kemp (1986), Rayner and Best (1995) -these three tests are not consistent against each fixed alternative- and, more recently, the tests in Novoa-Mu\~noz and Jim\'enez-Gamero (2014), and Novoa-Mu\~noz and Jim\'enez-Gamero (2016) (hereafter abbreviated to NJ (2014) and NJ (2016), respectively).
	
	The tests in NJ (2014) and NJ (2016) are consistent against each fixed alternative. The results in Janssen (2000) assert that the global power function of any nonparametric test is flat on balls of alternatives  except for alternatives coming from a finite-dimensional subspace. Therefore, it is interesting to propose new gof tests able to detect different sets of alternatives.
	
	The present work proposes a new consistent gof test for the BPD. To derive it we first show that the probability generating function (pgf) of the BPD is the only pgf that satisfies a certain system of partial differential equations. Therefore, under the null hypothesis, the empirical probability generating function (epgf), which is a consistent estimator of the pgf (see, e.g. NJ, 2014), should approximately satisfy such system. The proposed test statistic can be seen as a bivariate extension of the one in Baringhaus and Henze (1992) designed for testing gof to the univariate Poisson distribution.
	
	The asymptotic behavior of the proposed test under alternatives is shared with the ones in NJ (2014) or NJ (2016). An advantage of the test proposed in this paper over those in NJ (2014) and NJ (2016) is its speed for the delivery of results.

	In order to consistently approximate the null distribution of the test statistic, we propose to use a parametric bootstrap estimator. The finite-sample size performance of the test is numerically evaluated through a simulation study. The power of the test is compared with the tests mentioned above. There is no test yielding the highest power against each considered alternative, as expected from the results in Janssen (2000). In most cases, the power of the proposed test is quite close to the highest one; in other cases, the proposed test is the most powerful. In addition, from a computational point of view, the test proposed in this paper is more efficient than its competitors.

	Next we show the notation used in this work: all vectors are row vectors and $v^\top$ is the transpose of the row vector $v$; for any vector $v,\, v_k$ denotes its $k$th coordinate, and $\|v\|$ its Euclidean norm. We put $\,\mathbb{N}_0=\{0,1,2,3,\ldots\}$ and write $I_A\,$ for the indicator function of the set $A$; $P_{\theta}$ denotes the probability law of the BPD with parameter $\theta$; $P$ denotes the probability law of the data; $E_{\theta}$  denotes expectation regarding the probability function $P_{\theta}$; $E$  denotes expectation  with respect to the true probability function of the data; $P_*$ denote the  probability law, given the data; all limits in this work are taken as $n \rightarrow \infty; \mathop{\longrightarrow} \limits^{L}\ $ denotes convergence in distribution; $\mathop{\longrightarrow} \limits^{a.s.}\ $ denotes  almost sure (a.s.) convergence. For any function $h:S\subset\mathbb{R}^m \to \mathbb{R}$, for some fixed $m\in\mathbb{N}$, we will denote \[D^{k_1\cdots k_m}h(u)=\frac{\partial^k }{\partial u_1^{k_1}\cdots \partial u_m^{k_m}}\ h(u),\]  for each choice of nonnegative integers $k_1,\ldots,k_m$ such that $k=k_1+\cdots+k_m$.
	
	\section{Review of the existing tests for the BPD, and their multivariate extension}
	
	\subsection{Tests for the BPD}
	
	Let $\boldsymbol{X}_1=(X_ {11}, X_ {12}), \boldsymbol{X}_2=(X_ {21}, X_ {22}), \ldots, \boldsymbol{X}_n=(X_ {n1}, X_{n2})$ be independent identically distributed  (iid)  from  a random vector $\boldsymbol{X}=(X_1,X_2)$ taking values in $\mathbb{N}_0^2$. Based on the sample $\boldsymbol{X}_1,\boldsymbol{X}_2,\ldots,\boldsymbol{X}_n$, the objective is to test the hypothesis
	\[H_0 :\,\, (X_1,X_2)\sim BP(\theta_1,\theta_2,\theta_3), \ \text{for some}\ (\theta_1,\theta_2,\theta_3)\,\in\Theta,\]
	against the alternative
	\[ H_1 :\,\, (X_1,X_2)\nsim BP(\theta_1,\theta_2,\theta_3), \ \forall\,(\theta_1,\theta_2,\theta_3)\,\in\Theta,\]
	where $\Theta=\left\{\left(\theta_1,\theta_2,\theta_3\right)\in\mathbb{R}^3\ :\ \theta_1>\theta_3, \ \theta_2>\theta_3,\ \theta_3> 0\right\}$. From NJ (2014) the distribution of $\boldsymbol{X}=(X_1,X_2)$ is determined by its pgf $g(u)=E\left(u_1^{X_1}u_2^{X_2}\right),\, u=(u_1,u_2)\in [0,1]^2$, and the joint pgf of a random vector $\boldsymbol{X}\sim BP(\theta)$ is
	\begin{equation}\label{pgf_BPD}
	g(u;\theta)=E_{\theta}(u_1^{X_1}u_2^{X_2})=
	\exp\bigl\{\theta_1(u_1-1)+ \theta_2(u_2-1)+\theta_3(u_1-1)(u_2-1)\bigr\}.
	\end{equation}
	The empirical counterpart of pgf is epgf  of the data given by $g_n(u)=\frac{1}{n}\sum_{i=1}^n u_1^{X_{i1}}u_2^{X_{i2}}$.\\

	Next we will briefly expose three non-consistent tests that we found in the statistical literature, where $\bar{X}_1, \bar{X}_2, S^2_{X_1}$ and $S^2_{X_2}$ are the sample means and variances, respectively, $S^2_{X_1X_2}$ is the sample covariance, $r$ is the sample correlation coefficient and $\chi^2_{k, \alpha}$, for $0<\alpha<1$ and $k\in \mathbb{N}$, denotes the upper $\alpha$-percentile of the $\chi^2$ distribution with $k$ degrees of freedom.
	
	\subsubsection{Crockett test $\boldsymbol{T}$}
	The statistic $T$ (say) proposed by Crockett (1979) is based on a quadratic form in  $Z_{X_1}=S^2_{X_1}-\bar{X}_1$ and $Z_{X_2}=S^2_{X_2}-\bar{X}_2$. He shows that, under $H_0$, $T=ZV^{-1}Z^{\top}\ \mathop{\longrightarrow}\limits^{L}\ Y\sim \chi^2_{2}$, where $V$ denotes the matrix of variances and covariances of  $Z=(Z_{X_1}, Z_{X_2})$. Thus, the statistic and its critical region are given by
	\[ T= \frac{n}{2}\frac{\bar{X}_2^{2} \left(S^2_{X_1}-\bar{X}_1\right)^2- 2S^2_{X_1X_2}\left(S^2_{X_1}-\bar{X}_1\right) \left(S^2_{X_2}-\bar{X}_2\right)+ \bar{X}_1^{2} \left(S^2_{X_2}-\bar{X}_2\right)^2} {\bar{X}_1^{2} \bar{X}_2^{2}-S^4_{X_1X_2}}, \quad T \geq \chi^2_{2,\alpha}.\]

	\subsubsection{Test $\boldsymbol{I_B}$ of Loukas and Kemp}
	Loukas and Kemp (1986) developed a test based on what they called the bivariate dispersion index, $I_B=\frac{1}{1-\rho^2}\sum_{i=1}^n\left(W_{i1}^2-2\rho W_{i1}W_{i2}+W_{i2}^2\right),$
	where $W_{ik}=\frac{X_{ik}-\theta_k}{\sqrt{\theta_k}}, k=1,2, \, i=1,2,\ldots,n$ and $\ \rho=\frac{\theta_3}{\sqrt{\theta_1\theta_2}}$. If $\theta_1,\theta_2$ and $\theta_3$ are known, these authors show that $I_B$ is distributed approximately as a variable $\chi^2_{2n}$. If $\theta_1,\theta_2$ and $\theta_3$ are unknown, the statistic and its critical region are given by
	\[I_B=\frac{n(\bar{X}_2S^2_{X_1}-2S^2_{X_1X_2}+\bar{X}_1S^2_{X_2})} {\bar{X}_1\bar{X}_2-S^2_{X_1X_2}}, \quad I_B \geq \chi^2_{2n-3,\alpha}.\]
	
	\subsubsection{Test $\boldsymbol{NI_B}$ of Rayner and Best}
	Rayner and Best (1995) expressed the statistics of Loukas and Kemp (1986) as 
	$I_B=\frac{n}{1-\widehat{\rho}^{\,2}}\!\left(\frac{S_{X_1}^2}{\bar{X}_1}-2\, \frac{S_{X_1X_2}^2}{\bar{X}_1\bar{X}_2}+\frac{S_{X_2}^2} {\bar{X}_2}\right),$
	where $\widehat{\rho}=\frac{S_{X_1X_2}}{\sqrt{\bar{X}_1\bar{X}_2}}$ is an estimator of $\rho$. If $\widehat{\rho}^{\,2} >\frac{1}{2} \left(\frac{S_{X_1}^2}{\bar{X}_1}+\frac{S_{X_2}^2} {\bar{X}_2}\right)$, then $I_B<0$, and its distribution is not well approximated by a $\chi^2$. They proposed the statistic $$NI_B=\frac{n}{1-r^2}\left(\frac{S_{X_1}^2}{\bar{X}_1}-2\,r^2 \sqrt{\frac{S_{X_1}^2S_{X_2}^2}{\bar{X}_1\bar{X}_2}}+\frac{S_{X_2}^2} {\bar{X}_2}\right).$$ Under $H_0$, $NI_B$ is approximately distributed as $\chi^2_{2n-3}$ if $n$ is large. Therefore, $H_0$ is rejected if $NI_B \geq \chi^2_{2n-3,\alpha}.$\\

	Note that the statistical tests $T, I_B$ and $NI_B$ are not consistent, because they are based on the moments, specifically based on the fact that the first two population moments are equal. In contrast, the tests presented below are consistent.

	\subsubsection{Test $\boldsymbol{R_{n,w}}$ of Novoa-Mu\~noz and Jim\'enez-Gamero (2014)}
	From NJ (2014) the distribution of  $\boldsymbol{X}$ is uniquely determined by its pgf, $g(u)$, $u\in [0,1]^2$, a reasonable test for testing $H_0$ should reject the null hypothesis for large values of $$R_{n,w}=\int_0^1\int_0^1 G^2_n(u;\widehat{\theta}_n )w(u) du,$$ where 
	$G_n(u;\vartheta)= \sqrt{n}\{g_n(u)-g(u;\vartheta) \}$, $w(u)=u_1^{a_1}u_2^{a_2}$ is a measurable weight function
	$\forall u\in [0,1]^2$, $a_1, a_2\in (-1,\infty)$, and $\widehat{\theta}_n=\widehat{\theta}_n(\boldsymbol{X}_1, \ldots,\boldsymbol{X}_n)=(\widehat{\theta}_{1n}, \widehat{\theta}_{2n},\widehat{\theta}_{3n})$ is a consistent estimator of $\theta$.

	\subsubsection{Test $\boldsymbol{S_{n,w}}$ of Novoa-Mu\~noz and Jim\'enez-Gamero (2014)}
	Since the pgf $g(x)$ of the univariate Poisson distribution, with parameter $\lambda$, is the only pgf satisfying the differential equation $g'(x)=\lambda g(x)$, Baringhaus and Henze (1992) proposed a test statistic which is based on an empirical counterpart of this equation. With the aim of extending this result to the bivariate case, NJ (2014) proposed to reject $H_0$ for large values of $$S_{n,w}=n\int_0^1\int_0^1 \left\{B^2_{1n}(u;\widehat{\theta}_n )+B^2_{2n}(u;\widehat{\theta}_{n} )\right\}w(u)\ du,$$ where $w(u)=u_1^{a_1}u_2^{a_2}$, $a_1, a_2\in (-1,\infty)$, $\widehat{\theta}_n=(\widehat{\theta}_{1n}, \widehat{\theta}_{2n},\widehat{\theta}_{3n})$ is a consistent estimator of $\theta$, and
	\[
	\begin{array}{ll}
	B_{1n}(u;\widehat{\theta}_n )&=\displaystyle\frac{\partial }{\partial u_1}g_n(u_1,u_2)-\left\{\widehat{\theta}_{1n}+\widehat{\theta}_{3n}(u_2-1)\right\} g_n(u_1,u_2), \\[.3 cm]
	B_{2n}(u;\widehat{\theta}_{n} )&=\displaystyle\frac{\partial }{\partial u_2}g_n(u_1,u_2)-\left\{\widehat{\theta}_{2n}+\widehat{\theta}_{3n}(u_1-1)\right\} g_n(u_1,u_2),
	\end{array}
	\]
	should be close to 0 when $H_0$ is true. These functions are the empirical counterpart of the system of partial differential equations of Proposition 2 in NJ (2014).
	
	\subsubsection{Test $\boldsymbol{W_{n}}$ of Novoa-Mu\~noz and Jim\'enez-Gamero (2016)}
	When $H_0$ is true, NJ (2016) presented another interpretation of the fact that
	$S_{n,w}=n\int_0^1\int_0^1 \{B^2_{1n}(u;\widehat{\theta}_n )+B^2_{2n}(u;\widehat{\theta}_{n})\}w(u)\ du\approx 0$. Reasoning as Nakamura and P\'erez-Abreu (1993) for the univariate case and noting that $B_{kn}(u;\widehat{\theta}_n)=\sum_{r_1\geq 0}\sum_{r_2\geq 0}  b_k(r_1,r_2;\widehat{\theta}_n) u_1^{r_1} u_2^{r_2}$, $k=1,2$. NJ (2016) proposed to reject $H_0$ for large values of
	\[W_n=\sum_{r_1\geq 0}\sum_{r_2\geq 0}\{b^2_1(r_1,r_2;\widehat{\theta}_n)+b^2_2(r_1,r_2;\widehat{\theta}_n)\}=
	\sum_{r_1,r_2 = 0}^M\{b^2_1(r_1,r_2;\widehat{\theta}_n)+b^2_2(r_1,r_2;\widehat{\theta}_n)\},\]
	where $M=\max \{X_{(n)1}, X_{(n)2}\}$, $X_{(n)k}=\max_{1\leq i\leq n}X_{ik}$, $k=1,2$,
	\begin{align}
	b_1(r_1,r_2;\widehat{\theta}_n)&=(r_1+1)p_n(r_1+1,r_2)-(\widehat{\theta}_{1n}-\widehat{\theta}_{3n}) p_n(r_1,r_2)-\widehat{\theta}_{3n} p_n(r_1,r_2-1),\notag\\[.2 cm]
	b_2(r_1,r_2;\widehat{\theta}_n)&=(r_2+1)p_n(r_1,r_2+1)-(\widehat{\theta}_{2n}-\widehat{\theta}_{3n}) p_n(r_1,r_2)-\widehat{\theta}_{3n} p_n(r_1-1,r_2),\notag
	\end{align}
	and $p_n(r_1,r_2)=\frac{1}{n}\sum_{i=1}^n I_{\{X_{i1}=r_1,X_{i2}=r_2\}}$ is the relative frequency of the pair $(r_1,r_2)$, 
	
	\subsection{The general $m-$variate case}\label{The general m-variate case}
	For the multivariate case, for each integer $m>2$, let
	\[X_1=Y_1+Y_{m+1}, \,\, X_2=Y_2+Y_{m+1},\,\,  \ldots,  \,\, X_m=Y_m+Y_{m+1},\]
	where $Y_1, Y_2, \ldots, Y_{m+1}$ are independent Poisson random variables
	with means $ \theta'_1=\theta_1-\theta_{m+1}>0, \ldots, \theta'_m=\theta_m-\theta_{m+1}>0$ and  $\theta_{m+1}>0$, respectively. The joint distribution of the vector $(X_1,  X_2, \ldots, X_m)$ is called  a $m$-variate Poisson distribution with parameter $\theta=(\theta_1,\theta_2,\ldots, \theta_m,\theta_{m+1})$ (see Johnson, Kotz and Balakrishnan, 1997). The joint pgf of
	$(X_1,  X_2, \ldots, X_m)$ is
	\begin{equation}\label{fgp-multivariante}
	g(u;\theta)=
	\exp\left\{\sum_{i=1}^m \theta_i\left(u_i-1\right)+ \theta_{m+1}
	\left(\prod_{i=1}^m u_i-\sum_{i=1}^m u_i+m-1\right)\right\},\ \, \forall u \in \mathbb{R}^m.
	\end{equation}
	The empirical counterpart of pgf is epgf of the data given by \begin{equation}\label{efgp-multivariante}
	g_n(u)=\frac{1}{n}\sum_{i=1}^n u_1^{X_{i1}}\cdots u_m^{X_{im}},\ \, u=(u_1,\ldots,u_m).
	\end{equation}
	Now, the objective is to test the hypothesis
	\[H_{0m}:(X_1,  X_2, \ldots, X_m)\mbox{  has a $d$-variate Poisson distribution}.\]\\
	The tests proposed by Crockett (1979), Loukas and Kemp (1986), and Rayner and Best (1995) do not have a multivariate extension. However, NJ (2014) and NJ (2016) proposed a natural extension of their tests, which will be presented below.
	
	\subsubsection{Test $\boldsymbol{R_{m,n,w}}$ of Novoa-Mu\~noz and Jim\'enez-Gamero (2014)}
	NJ (2014) affirmed that the extension of the test $R_{n,w}$ is direct, it is enough to consider pgf $g(u;\theta)$ as in (\ref{fgp-multivariante}), epgf $g_n(u;\widehat{\theta}_n)$ as in (\ref{efgp-multivariante}) and $w(u)$ is a measurable nonnegative weight function with finite integral over $[0,1]^m$.

	\subsubsection{Test $\boldsymbol{S_{m,n,w}}$ of Novoa-Mu\~noz and Jim\'enez-Gamero (2014)}
	To test $H_{0m}$, NJ (2014) considered the test statistic
	\[S_{m,n,w}=n\int_{[0,1]^m} \left\{B^2_{1n}(u;\widehat{\theta}_n )+\cdots+B^2_{mn}(u;\widehat{\theta}_{n} )\right\}w(u)\ du,\]
	where $w(u)$ is a measurable nonnegative weight function with finite integral over $[0,1]^m$, and
	\[
	B_{in}(u;\widehat{\theta}_n)=\displaystyle\frac{\partial }{\partial u_i}g_n(u)-\left\{\widehat{\theta}_{i,n}+\widehat{\theta}_{m+1,n}\left(\prod _{j\neq i }u_j-1\right)\right\} g_n(u),\ 1\leq i \leq m,\ g_n(u) \text{ as in (\ref{efgp-multivariante})}.\]

	\subsubsection{Test $\boldsymbol{W_{m,n}}$ of Novoa-Mu\~noz and Jim\'enez-Gamero (2016)}
	With the aim of extending $W_{n}$ to the multivariate case, NJ (2016) proposed the following statistic for testing $H_{0m}$,
	\[W_{m,n}=\sum_{r_1,r_2,\ldots,r_m\geq 0}\left\{\sum_{j=1}^m b^2_j(r_1,r_2,\ldots,r_m;\widehat{\theta}_n)\right\}=
	\sum_{r_1,r_2,\ldots,r_m= 0}^M\left\{\sum_{j=1}^m b^2_j(r_1,r_2,\ldots,r_m;\widehat{\theta}_n)\right\},\]
	where $M=\max \{X_{(n)1}, X_{(n)2},\ldots,X_{(n)m}\}$, $X_{(n)k}=\max_{1\leq i\leq n}X_{ik}$, $1\leq k\leq m$, and 
	\begin{align}
	b_j(r_1,\ldots,r_m;\widehat{\theta}_n)&=(r_j+1) p_n(r_1,\ldots,r_{j-1},r_j+1,r_{j+1},\ldots,r_m)-(\widehat{\theta}_{jn}- \widehat{\theta}_{m+1,n}) p_n(r_1,\ldots,r_m)\notag \\[.2 cm]
	&\hspace{4mm}-\widehat{\theta}_{m+1,n} \, p_n(r_1-1,\ldots,r_{j-1}-1,r_j,r_{j+1}-1,\ldots,r_m-1),\ \, 1\leq j\leq m,\notag
	\end{align}
	and $p_n(r_1,\ldots,r_m)=\frac{1}{n}\sum_{i=1}^n I_{\{X_{i1}=r_1,\ldots,X_{im}=r_m\}}$
	is the relative frequency of $(r_1,\ldots,r_m)$.

	\section{A new characterization of the BPD}\label{A characterization of the BPD}
	In order to obtain a new test to test the hypothesis $H_0$ against the alternative $H_1$ and based on the fact that the distribution of $\boldsymbol{X}=(X_1,X_2)$ is determined by its pgf, we give a different characterization for the BPD.
	\begin{proposition}\label{Soluc-EDP-conjunta-y-marginales} Let $g(u_1,u_2;\theta)\,$ be as defined in (\ref{pgf_BPD}). Then $g(u_1,u_2;\theta)\,$ is the only pgf satisfying the following system of partial differential equations
		\begin{equation}\label{EDP-conjunta-y-marginales}
		\left.
		\begin{array}{rcl}
		D_1(u;\theta) &=&\displaystyle\frac{\partial g(u_1,1)}{\partial u_1}-\theta_1\, g(u_1,1)=0,\\[.3 cm]
		D_2(u;\theta)&=&\displaystyle\frac{\partial g(1,u_2)}{\partial u_2}-\theta_2\, g(1,u_2)=0,\\[.3 cm]
		D_3(u;\theta)&=&\displaystyle\frac{\partial^2 g(u_1,u_2)}{\partial u_1 \,\partial u_2}-f(u_1,u_2;\theta)\, g(u_1,u_2)=0,  \end{array}
		\right\}
		\end{equation}
		where $\,f(u_1,u_2;\theta)=\theta_3+ \{\theta_2+\theta_3(u_1-1)\}\{\theta_1+\theta_3(u_2-1)\}$.
	\end{proposition}
	
	The system of equations (\ref{EDP-conjunta-y-marginales}) has the following nice interpretation: first and the second equation characterize the marginal distributions, i. e., they are equivalent to saying that the marginal distributions are univariate Poisson; the last equation characterizes the dependence structure.

	By Proposition 1 in NJ (2014), $g(u)$ and its derivatives can be consistently estimated by the epgf and the derivatives of the epgf, respectively. Thus, if $H_0$ is true, then the functions
	$$  \begin{array}{rcl}
	D_{1n}\bigl(u;\widehat{\theta}_{n} \bigr)&=&\displaystyle\frac{\partial g_n(u_1,1)}{\partial u_1}-\widehat{\theta}_{1n}\ g_n(u_1,1), \\[.3 cm]
	D_{2n}\bigl(u;\widehat{\theta}_{n} \bigr)&=&\displaystyle\frac{\partial g_n(1,u_2)}{\partial u_2}-\widehat{\theta}_{2n}\ g_n(1,u_2),\\[.3 cm]
	D_{3n}\bigl(u;\widehat{\theta}_n \bigr)&=&\displaystyle\frac{\partial^2 g_n(u_1,u_2)}{\partial u_1 \,\partial u_2}-f(u_1,u_2;\widehat{\theta}_n)\, g_n(u_1,u_2),
	\end{array}$$
	should be close to $0$, $\forall (u_1,u_2) \in [0,1]^2$, where $\widehat{\theta}_n=(\widehat{\theta}_{1n}, \widehat{\theta}_{2n}, \widehat{\theta}_{3n})$ is a consistent estimator of $\theta$ and
	$g_n(u_1,u_2)$ is the epgf associated with the data, i. e.,
	\[g_n(u_1,u_2)=\frac{1}{n}\sum_{i=1}^n u_1^{X_{i1}}u_2^{X_{i2}}.\]
	
	Thus, to test $H_0$ we consider the following test statistic
	\begin{equation}\notag
	T_{n,w}=n\int_0^1\int_0^1 \left\{D^2_{1n}\bigl(u;\widehat{\theta}_n \bigr)+D^2_{2n}\bigl(u;\widehat{\theta}_{n} \bigr)+D^2_{3n}\bigl(u;\widehat{\theta}_{n} \bigr)\right\}w(u)\ du,
	\end{equation}
	where $w(u)$ is a non-negative function on
	$[0,1]^2$.
	
	In order to give a sound justification of  $\,T_{n,w}$ as a test statistic for testing $H_0$ we next derive its almost sure limit.

	\begin{theorem}\label{Convergencia-de-T_n,w(theta-hat)}
		Let $\boldsymbol{X}_1, \boldsymbol{X}_2,\ldots,\boldsymbol{X}_n$ be iid from  $\boldsymbol{X}=(X,Y)\in\mathbb{N}_0^2$  with pgf $g(u)$ such that $\frac{\partial g(u_1,1)}{\partial u_1}$, $\frac{\partial g(1,u_2)}{\partial u_2}$ and $\frac{\partial^2 g(u_1,u_2)}{\partial u_1 \,\partial u_2}$, exist and are continuous functions on a region containing $[0,1]^2$. If
		$\widehat{\theta}_n\mathop{\longrightarrow}\limits^{a.s.} \theta$, for some $\,\theta\in \mathbb{R}^3$, then
		$$\frac{T_{n,w}}{n}\ \mathop{\longrightarrow}\limits^{\!a.s.}\int_0^1\!\int_0^1 \left\{D^2_1(u;\theta)+D^2_2(u;\theta)+D^2_3(u;\theta)\right\}w(u)\, du=\eta(g;\theta)\geq 0.$$
	\end{theorem}
	
	Note that if $w>0$ almost everywhere (a.e.) on $[0,1]^2$, then $\eta(g;\theta)=0$ if and only if $H_0$ is true. Therefore, a reasonable test for testing $H_0$ should reject the null hypothesis for large values of $T_{n,w}$. Now, to determine what are large values of $T_{n,w}$, we must calculate its null distribution, or at least an approximation to it. Clearly, the null distribution of $T_{n,w}$ is unknown. A classical way of approximating the null distribution of a test statistic is through its asymptotic null distribution. The next section studies this issue.
	
	\section{A bootstrap estimator of the null distribution}\label{Null bootstrap distribution}
	In order to derive the asymptotic null distribution of the test statistic $T_{n,w}$ we will assume that the estimator $\widehat{\theta}_n$ is asymptotically linear, as expressed in Assumption 1 in NJ (2014) and we will consider the separable Hilbert space 
	$$\mathcal{H}=\{\varphi:[0,1]^2\rightarrow\mathbb{R},\ \ \text{with}\ \ \|\varphi\|_{_\mathcal{H}}^{\,2}=\int_0^1\int_0^1 \varphi^2(u)\,w(u)\,du<\infty\}.$$
	
	In this framework, $T_{n,w}$ can be expressed as   $T_{n,w}= \|Z_{1n}\|_{_\mathcal{H}}^{\,2} +\|Z_{2n}\|_{_\mathcal{H}}^{\,2} +\|Z_{3n}\|_{_\mathcal{H}}^{\,2}$, with
	\[Z_{kn}(u)=\frac{1}{\sqrt{n}}\sum_{i=1}^n R_k\bigl(\boldsymbol{X}_i;\widehat{\theta}_n; u\bigr), \ \ k=1,2,3,\]
	where, for $ 1\leq i\leq n$,
	\begin{align}
	R_1\bigl(\boldsymbol{X}_i;\widehat{\theta}_n; u\bigr)&= X_{i1}\,I_{\{X_{i1}\geq 1\}}\,u_1^{X_{i1}-1}-\widehat{\theta}_{1n}\,u_1^{X_{i1}}, \notag\\[.1 cm]
	R_2\bigl(\boldsymbol{X}_i;\widehat{\theta}_n; u\bigr)&= X_{i2}\,I_{\{X_{i2}\geq 1\}}\,u_2^{X_{i2}-1}-\widehat{\theta}_{2n}\,u_2^{X_{i2}}, \notag\\[.1 cm]
	R_3\bigl(\boldsymbol{X}_i;\widehat{\theta}_n; u\bigr)&= X_{i1}\,X_{i2}\,I_{\{X_{i1}X_{i2}\geq 1\}}\,u_1^{X_{i1}-1}\, u_2^{X_{i2}-1}-f\bigl(u;\widehat{\theta}_n\bigr)u_1^{X_{i1}}\, u_2^{X_{i2}},\notag
	\end{align}
	The next result gives the asymptotic null distribution of $T_{n,w}$.
	
	\begin{theorem}\label{TeoConvDebil-T_n,w}
		Let  $\,\boldsymbol{X}_{1},\boldsymbol{X}_{2},\ldots, \boldsymbol{X}_{n}\,$ be iid from $\,\boldsymbol{X}=(X_{1},X_{2})\sim BP(\theta)$. Suppose that Assumption 1 in NJ (2014) holds and that $\widehat{\theta}_n\mathop{\longrightarrow}\limits^{a.s.} \theta$. Then
		\[T_{n,w}= \|W_{1n}\|_{_\mathcal{H}}^{\,2}+ \|W_{2n}\|_{_\mathcal{H}}^{\,2}+\|W_{3n}\|_{_\mathcal{H}}^{\,2}+r_n,\]
		where $\ P_{\theta}(|r_n|>\varepsilon)\to 0$, $\forall \varepsilon>0$,
		\[W_{kn}(u)=\frac{1}{\sqrt{n}}\sum_{i=1}^n W^0_k(\boldsymbol{X}_i;\theta; u),\ \ k=1,2,3,\]
		\begin{align}
		W^0_1(\boldsymbol{X}_i;\theta; u)&=X_{i1}\,I_{\{X_{i1}\geq 1\}}\,u_1^{X_{i1}-1}-\theta_{1}\,u_1^{X_{i1}}- g(u_1,1;\theta)\,\boldsymbol{\ell}\left(\boldsymbol{X}_{i}; \theta\right)(1,0,0)^\top, \notag \\[.25 cm]
		W^0_2(\boldsymbol{X}_i;\theta; u)&=X_{i2}\,I_{\{X_{i2}\geq 1\}}\,u_2^{X_{i2}-1}-\theta_{2}\,u_2^{X_{i2}}- g(1,u_2;\theta)\,\boldsymbol{\ell}\left(\boldsymbol{X}_{i}; \theta\right)(0,1,0)^\top, \notag \\[.25 cm]
		W^0_3(\boldsymbol{X}_i;\theta; u)&=X_{i1}X_{i2}I_{\{X_{i1}X_{i2}\geq 1\}}u_1^{X_{i1}\!-1}u_2^{X_{i2}\!-1}-f(u;\theta) u_1^{X_{i1}}u_2^{X_{i2}}-g(u;\theta) \boldsymbol{\ell}\left(\boldsymbol{X}_{i}; \theta\right)B^\top(u;\theta), \notag
		\end{align}
		$1\leq i \leq n$, $B(u;\theta)=\!\left(b_1(u;\theta), b_2(u;\theta),b_3(u;\theta)\right),$ where $\,b_1(u;\theta)=\theta_{2}+\theta_{3}(u_1-1)$, $b_2(u;\theta)=\theta_{1}+\theta_{3}(u_2-1)\,$ and $\ b_3(u;\theta)=1+\theta_{1}(u_1-1)+\theta_{2}(u_2-1)+2\theta_{3}(u_1-1)(u_2-1)$.
		Moreover,
		\[T_{n,w}\ \mathop{\longrightarrow} \limits^{\!L}\ \sum_{j\geq 1}\lambda_j\,\chi^2_{1j}\,,\]
		where $\chi^2_{11},\chi^2_{12},\ldots$ are independent $\chi^2$ variates with one degree of freedom and the set $\{\lambda_j\}$ are the non-null eigenvalues of the operator $C(\theta)$ defined on the function space $\{\tau:\mathbb{N}_0^2\to \mathbb{R}, \ \text{such that} \ E_{\theta}\!\left\{\tau^2(\boldsymbol{X})\right\}<\infty,\forall \theta\in\Theta\}$, as follows
		\begin{equation}\label{expresion_del_operador_C(theta)}
		C(\theta)\, \tau(\boldsymbol{x})= E_{\theta}\{h(\boldsymbol{x},\boldsymbol{Y};\theta)\, \tau(\boldsymbol{Y})\},
		\end{equation}
		\noindent with \begin{equation}\label{nucleo_de_Tnw_aproximado}
		h(\boldsymbol{x},\boldsymbol{y};\theta)=\int_0^1\int_0^1 \sum_{k=1}^3 W_k^0(\boldsymbol{x}; \theta; u)\, W_k^0(\boldsymbol{y}; \theta; u)\,w(u)\,du.
		\end{equation}
	\end{theorem}
	
	The asymptotic null distribution of $T_{n,w}$ does not provide a useful approximation to its null distribution since it depends on the unknown true value of $\theta$. This could be overcome by replacing $\theta$ by $\widehat{\theta}_n$. But the greatest difficulty is to determine the set $\{\lambda_j\}$, since, in general, calculating the eigenvalues of an operator is not an easy task and in our case we must also obtain expression (\ref{nucleo_de_Tnw_aproximado}), which is not easy to derive. So, we next consider another way of approximating the  null distribution of the test statistic, the bootstrap.\\

	The following result proves that the bootstrap method consistently approximates the null distribution of $T_{n,w}$, for which we require the Assumption 2 in NJ (2014) and the previous explanations for that assumption.
	\begin{theorem}\label{Teo_Cons_Boot_T_n,w}
		Let $\,\boldsymbol{X}_{1},\ldots, \boldsymbol{X}_{n}\,$ be iid random vectors from $\boldsymbol{X}=(X_{1},X_{2}) \in \mathbb{N}_0^2$. Suppose that  Assumption 2 in NJ (2014) holds, $\widehat{\theta}_n\mathop{\longrightarrow}\limits^{a.s.} \theta$, for some $\theta \in \Theta$. Then
		$$\sup_{x\in\,\mathbb{R}} \left|P_*\!\left(T^*_{n,w} \leq x\right)-P_{\theta}\!\left(T_{n,w} \leq x\right)\right|\ \mathop{\longrightarrow}\limits^{a.s.}\ 0.$$
	\end{theorem}
	
	It is important to note that analogous comments follow those given after Theorem 2 in NJ (2014) and the test function for our case is presented below.
	
	Let
	$t^*_{n,w,\alpha}=\inf\{x:P_*\bigl(T^*_{n,w}\geq x\bigr)\leq \alpha\}.$
	From Theorem \ref{Teo_Cons_Boot_T_n,w}, the test function
	$$\Psi^*=\left\{
	\begin{array}{ll}
	1, & \text{if}\ T_{n,w}\geq t^*_{n,w,\alpha}\,, \\[.2 cm]
	0, & \text{otherwise},
	\end{array}
	\right.$$
	or equivalently, the test that rejects $H_0$ when
	$p^*=P_*\!\left(T^*_{n,w} \geq T_{obs}\right)\leq \alpha$, is asymptotically correct, in the sense that the type I error is asymptotically equal to the nominal value $\alpha$, where $T_{obs}$ is the observed value of the test statistic $T_{n,w}$.

	\section{Behaviour against alternatives}\label{Alternatives}
	As an immediate consequence of Theorems \ref{Convergencia-de-T_n,w(theta-hat)}, \ref{TeoConvDebil-T_n,w} and \ref{Teo_Cons_Boot_T_n,w}, the next result gives the asymptotic power of the test $\Psi^*$ against fixed alternatives.
	
	\begin{corollary}\label{Altern-fijas}
		Let $\boldsymbol{X}_1, \boldsymbol{X}_2,\ldots,\boldsymbol{X}_n$ be iid from  $\boldsymbol{X}\in\mathbb{N}_0^2$ with pgf $g(u)$. Suppose that assumptions in Theorems \ref{Convergencia-de-T_n,w(theta-hat)} and \ref{Teo_Cons_Boot_T_n,w} hold. If $\eta(g;\theta)>0$, then
		$P\bigl(\Psi^*=1\bigr) \to 1.$
	\end{corollary}
	
	As commented after Theorem \ref{Convergencia-de-T_n,w(theta-hat)}, a simple way to ensure that $\eta(g;\theta)>0$, $\forall
	(X_1,X_2)\nsim BP(\theta)$, $ \forall\,(\theta)\in \Theta$, and thus the consistency against any fixed alternative,  is by choosing the weight function $w$ positive a.e. on $[0,1]^2$.\\
	
	For the local power, the next result ensures that the test $\Psi^*$ is able to detect alternatives as defined in (11) in NJ (2014), which converge to the BPD at the rate $n^{-1/2}$. With this aim, let $\{\phi_j\}$ be the set of orthonormal eigenfunctions corresponding to the eigenvalues $\{\lambda_j\}$ of the operator $C(\theta)$ given in (\ref{expresion_del_operador_C(theta)}).

	\begin{theorem}\label{Altern_contiguas_T}
		Let $\boldsymbol{X}_1, \boldsymbol{X}_2,\ldots,\boldsymbol{X}_n$ be iid from  $\boldsymbol{X} \in \mathbb{N}_0^2$, with pmf $P_n(x,y)$ as defined in (11) in NJ (2014). Suppose that Assumptions 1 and 3 in NJ (2014) hold. Then
		\[T_{n,w}\ \mathop{\longrightarrow} \limits^{\!L} \ \sum_{k=1}^{\infty} \lambda_k\left(Z_k+c_k\right)^2,\]
		\noindent where $\,c_k=\mathop{\sum}\limits_{x,\,y}b(x,y)\,\phi_k(x,y)\,$ and $Z_1,Z_2,\ldots\,$ are independent standard normal variates.
	\end{theorem}

	\section{Some computational issues}\label{issues}
	
	\subsection{On the calculation of the test statistic}
	
	Using the weight function (11) in NJ (2014) we obtained the following expression of our statistic.
	\[T_{n,w}=\frac{1}{n}\sum_{i=1}^n\sum_{j=1}^n \bigl(T_{ij}^1+T_{ij}^2+T_{ij}^3\bigr),\]
	\begin{align}
	T_{ij}^k&=\frac{1}{a_k+1}\left\{\!\frac{X_{ik}\,I_{\!B_{ik}}\,X_{jk}\, I_{\!B_{jk}}}{X_{ik}+X_{jk}+a_k-1}-\frac{\widehat{\theta}_{kn}\left(X_{ik} I_{\!B_{ik}}\!+\!X_{jk}I_{\!B_{jk}}\right)}{X_{ik}+X_{jk}+a_k}+\frac{\widehat{\theta}_{kn}^{\ 2}}{X_{ik}+X_{jk}+a_k+1}\right\}\!,\, k=1,2,\notag\\[.1 cm]
	T_{ij}^3&=\frac{X_{\!i1}\,I_{\!B_{i1}}\,X_{\!i2}\,I_{\!B_{i2}}\,X_{\!j1}\,I_{\!B_{j1}}\, X_{\!j2}\,I_{\!B_{j2}}}{(X_{i1}+X_{j1}+a_1-1)(X_{i2}+X_{j2}+a_2-1)} - \frac{2\left\{\!\bigl(\widehat{\theta}_{\!1n}\!- \! \widehat{\theta}_{3n}\bigr)\!\bigl(\widehat{\theta}_{2n}\!- \! \widehat{\theta}_{3n}\bigr)\!+\!\widehat{\theta}_{3n}\! \right\}\!X_{\!j1}I_{\!B_{j1}}\!X_{\!j2}I_{\!B_{j2}}} {(X_{i1}+X_{j1}+a_1)(X_{i2}+X_{j2}+a_2)} \notag\\[.1 cm]
	&\ \ -\frac{2\,\widehat{\theta}_{3n}\bigl(\widehat{\theta}_{2n}- \widehat{\theta}_{3n}\bigr)X_{j1}\,I_{\!B_{j1}}\,X_{j2}\,I_{\!B_{j2}}} {(X_{i1}+X_{j1}+a_1)(X_{i2}+ X_{j2}+a_2+1)}-\frac{2\,\widehat{\theta}_{3n}\bigl(\widehat{\theta}_{1n}- \widehat{\theta}_{3n}\bigr)X_{j1}\,I_{\!B_{j1}}\,X_{j2}\,I_{\!B_{j2}}} {(X_{i1}+X_{j1}+a_1+1)(X_{i2}+ X_{j2}+a_2)}\notag\\[.1 cm]
	&\ \ + \frac{\bigl\{\bigl(\widehat{\theta}_{1n}- \widehat{\theta}_{3n}\bigr)\bigl(\widehat{\theta}_{2n}- \widehat{\theta}_{3n}\bigr)\!+\widehat{\theta}_{3n}\bigr\}^2- 2\,\widehat{\theta}_{3n}^{\ 2} \, X_{j1}\,I_{\!B_{j1}}\,X_{j2}\,I_{\!B_{j2}}} {(X_{i1}+X_{j1}+a_1+1)(X_{i2}+X_{j2}+a_2+1)}\notag\\[.1 cm]
	&\ \ + \!\frac{2\widehat{\theta}_{3n}\! \left\{\bigl(\widehat{\theta}_{1n}\!\!- \widehat{\theta}_{3n}\!\bigr)\bigl(\widehat{\theta}_{2n}\!\!-\! \widehat{\theta}_{3n}\!\bigr)\!+\!\widehat{\theta}_{3n}\!\right\}\! \bigl(\widehat{\theta}_{2n}\!\!- \! \widehat{\theta}_{3n}\!\bigr)}{(X_{i1}\!+\!X_{j1}\!+a_1\!+\!1)(X_{i2}\!+\! X_{j2}\!+\!a_2\!+\!2)}+\frac{2\widehat{\theta}_{3n}\! \left\{\bigl(\widehat{\theta}_{1n}\!-\! \widehat{\theta}_{3n}\!\bigr)\bigl(\widehat{\theta}_{2n}\!-\! \widehat{\theta}_{3n}\!\bigr)\!+\!\widehat{\theta}_{3n}\!\right\}\! \bigl(\widehat{\theta}_{1n}\!\!- \! \widehat{\theta}_{3n}\!\bigr)\!}{(X_{i1}\!+\!X_{j1}\!+\!a_1\!+\!2)(X_{i2}\!+\! X_{j2}\!+\!a_2\!+\!1)}\notag\\[.1 cm]
	&\ \ +\! \frac{\widehat{\theta}_{3n}^{\ 2}\bigl(\widehat{\theta}_{2n}- \widehat{\theta}_{3n}\!\bigr)^{\!2}}{(X_{\!i1}\!+\!X_{\!j1}\!+a_1+1) (X_{i2}\!+\!X_{j2}\!+a_2+3)\!}+\!\frac{\widehat{\theta}_{3n}^{\ 2}\bigl(\widehat{\theta}_{1n}- \widehat{\theta}_{3n}\!\bigr)^{\!2}}{(X_{\!i1}\!+\!X_{\!j1}\!+a_1+3) (X_{i2}\!+\!X_{j2}\!+a_2+1)\!} \notag\\[.1 cm]
	&\ \ + \!\frac{2\,\widehat{\theta}_{3n}^{\ 2}\!\left\{2\bigl(\widehat{\theta}_{1n}-\widehat{\theta}_{3n} \!\bigr)\bigl(\widehat{\theta}_{2n}-\widehat{\theta}_{3n}\! \bigr)+\widehat{\theta}_{3n}\right\}}{(X_{\!i1}\!+\!X_{\!j1}\!+a_1+2) (X_{i2}\!+\!X_{j2}\!+a_2+2)\!}+\!\frac{\widehat{\theta}_{3n}^{\ 4}}{(X_{\!i1}\!+\!X_{\!j1}\!+a_1+3) (X_{i2}\!+\!X_{j2}\!+a_2+3)\!} \notag\\[.1 cm]
	&\ \ +\!\! \frac{2\,\widehat{\theta}_{3n}^{\ 3}\bigl(\widehat{\theta}_{2n}- \widehat{\theta}_{3n}\!\bigr)\!}{(X_{\!i1}\!+\!X_{\!j1}\!+a_1\!+2) (X_{i2}\!+\!X_{j2}\!+a_2\!+3)\!}+\frac{2\,\widehat{\theta}_{3n}^{\ 3}\bigl(\widehat{\theta}_{1n}- \widehat{\theta}_{3n}\!\bigr)\!}{(X_{\!i1}\!+\!X_{\!j1}\!+a_1\!+3) (X_{i2}\!+\!X_{j2}\!+a_2\!+2)\!}\,,\notag
	\end{align}
	where $\, B_{rs}=\{X_{rs}\geq 1\},\, \ 1\leq r \leq n, \ s=1,2$.

	\subsection{On the calculation of the null bootstrap distribution estimator}
	In practice, the exact bootstrap estimator of  the null distribution of  $T_{n,w}$ cannot be calculated, we will approximate it by simulation following the parametric bootstrap procedure (\textbf{PB algorithm}) given in section 4.1 in NJ (2016).
	
	\section{Numerical results}\label{Simulations}
	The properties studied so far describe the  behavior of the  proposed test  for very large samples.  We carried a simulation experiment in order to study the goodness of the bootstrap approximation  as well as to compare the power of the proposed test with other tests for finite sample sizes.
	We briefly describe it in this section and display a summary of the results obtained. All computations were performed by using programs written in the R language.
	
	\subsection{Simulated data}
	In addition to the test proposed in this paper, $T_{n,a}$, we also considered the tests  given in Crockett (1979) (denoted by $T$, see subsection 2.1.1), Loukas and Kemp (1986) (denoted by $I_B$, see subsection 2.1.2), Rayner and Best (1995) (denoted by $NI_B$, see subsection 2.1.3), NJ (2014) (denoted by $R_{n,a}$ and $S_{n,a}$, see subsections 2.1.4 and 2.1.5, respectively) and NJ (2016) (denoted by $W_{n}$, see subsection 2.1.6).
	
	We studied the goodness of the proposed bootstrap approximations to the null distribution of the test statistic for finite sample sizes. With this aim, we generated 1,000 samples of size $n=30(20)70$ from $BP(\theta_1,\theta_2,\theta_3)$, with $\theta_1=\theta_2=1$ and $\theta_3$ such that the  correlation coefficient,
	$\rho={\theta_3}/{\sqrt{\theta_1\,\theta_2}}$, equals 0.25, 0.5 and 0.75. To estimate $\theta$ we employed the maximum likelihood method. Then we approximate the $p$-values bootstrap of the proposed tests with weight function (11) in NJ (2014) for $a=(a_1,a_2)\in \{(0,0), (1,0)\}$ and 500 bootstrap samples, as well as the (asymptotic) $p$-values associated with the test statistics $T$, $I_B$ and $NI_B$.
	
	We repeated the above experiment for $\theta_1=1.5$, $\theta_2=1$ and $\theta_3$ such that the  correlation coefficient (approximately) equals 0.25, 0.5 and 0.75. In this case, since $\theta_1\neq \theta_2$, we considered  $(a_1,a_2)\in \{(0,0), (1,0), (0,1)\}$ for $R_{n,a}$, $S_{n,a}$ and $T_{n,a}$ in order to examine the effect of giving different weight to each component when they have different means.
	
	Tables \ref{type I error_1} and \ref{type I error_2} display the fraction of estimated $p$-values less than or equal to $0.05$ and $0.10$, which are the estimated type I error probabilities for $\alpha=0.05$ and $0.10$ (denoted as f05 and f10 in the tables), respectively.

\renewcommand{\thetable}{\Roman{table}}
\begin{table}\begin{center}
\caption{ Simulation results for the probability of type I error, $\theta_1=\theta_2=1$.}\label{type I error_1}\vspace{5pt}
{\footnotesize
\begin{tabular}{|c|l|cc|c|cc|c|cc|c|}
\cline{3-11}\multicolumn{2}{c|}{}
  &\multicolumn{3}{c|}{$n=30$}& \multicolumn{3}{c|}{$n=50$} & \multicolumn{3}{c|}{$n=70$}\\ 
  \hline 
  $\theta_3=\rho$ & Test & f05 & f10 & $KS$ & f05 & f10 & $KS$ & f05 & f10 & $KS$\\
  \hline
  0.25 &$R_{n,(0,0)}$ & 0.037 & 0.087 & 0.863178 & 0.047 & 0.100 & 0.257432 & 0.044 & 0.086 & 0.111356\\
  &$S_{n,(0,0)}$ & 0.046 & 0.089 & 0.934732 & 0.046 & 0.089 & 0.818621 & 0.045 & 0.092 & 0.508494\\
  &$T_{n,(0,0)}$ & 0.043 & 0.087 & 0.329116 & 0.041 & 0.089 & 0.508494 & 0.047 & 0.098 & 0.718379\\
  \cline{2-11}
  &$R_{n,(1,0)}$ & 0.035 & 0.099 & 0.329116 & 0.047 & 0.103 & 0.818621 & 0.038 & 0.084 & 0.329116\\
  &$S_{n,(1,0)}$ & 0.034 & 0.090 & 0.902243 & 0.046 & 0.097 & 0.960002 & 0.041 & 0.089 & 0.508494\\
  &$T_{n,(1,0)}$ & 0.042 & 0.094 & 0.329116 & 0.038 & 0.095 & 0.129364 & 0.046 & 0.094 & 0.369615\\
  \cline{2-11}
&$W_{n}$       & 0.022 & 0.056 & 1.00e-05 & 0.033 & 0.078 & 0.111356 & 0.038 & 0.090 & 0.612128\\
  \cline{2-11}
 & $T$           & 0.011 & 0.031 & $<$ 2.2e-16 & 0.046 & 0.092 & 0.060937 & 0.013 & 0.038 & $<$ 2.2e-16 \\
 & $I_B$         & 0.027 & 0.061 & $<$ 2.2e-16 & 0.098 & 0.144 & 0.001642 & 0.022 & 0.054 & $<$ 2.2e-16 \\
 & $NI_B$        & 0.010 & 0.034 & $<$ 2.2e-16 & 0.068 & 0.111 & 0.003452 & 0.013 & 0.033 & $<$ 2.2e-16 \\
  \hline
  \hline
 0.50 & $R_{n,(0,0)}$ & 0.048 & 0.112 & 0.129364 & 0.044 & 0.106 & 0.197933 & 0.045 & 0.098 & 0.559560\\
  &$S_{n,(0,0)}$ & 0.041 & 0.094 & 0.049545 & 0.049 & 0.099 & 0.257432 & 0.049 & 0.099 & 0.413150\\
  &$T_{n,(0,0)}$ & 0.041 & 0.098 & 0.257432 & 0.046 & 0.085 & 0.149677 & 0.055 & 0.111 & 0.413150\\
  \cline{2-11}
  &$R_{n,(1,0)}$ & 0.051 & 0.101 & 0.129364 & 0.044 & 0.097 & 0.863178 & 0.047 & 0.109 & 0.172476\\
  &$S_{n,(1,0)}$ & 0.042 & 0.099 & 0.069329 & 0.050 & 0.095 & 0.291736 & 0.046 & 0.099 & 0.197933\\
  &$T_{n,(1,0)}$ & 0.044 & 0.095 & 0.413150 & 0.043 & 0.084 & 0.129364 & 0.051 & 0.112 & 0.459543\\
  \cline{2-11}
 &$W_{n}$       & 0.022 & 0.061 & 0.013476 & 0.032 & 0.077 & 0.111356 & 0.037 & 0.081 & 0.111356\\
  \cline{2-11}
  &$T$           & 0.026 & 0.049 & 1.40e-06 & 0.024 & 0.039 & $<$ 2.2e-16 & 0.021 & 0.053 & 0.000179 \\
  &$I_B$         & 0.088 & 0.125 & $<$ 2.2e-16 & 0.073 & 0.119 & $<$ 2.2e-16 & 0.051 & 0.081 & $<$ 2.2e-16 \\
  &$NI_B$        & 0.036 & 0.074 & 5.00e-07 & 0.018 & 0.049 & $<$ 2.2e-16 & 0.007 & 0.035 & $<$ 2.2e-16 \\
  \hline
  \hline
 0.75 &$R_{n,(0,0)}$ & 0.043 & 0.089 & 0.718379 & 0.060 & 0.112 & 0.902243 & 0.050 & 0.114 & 0.508494\\
  &$S_{n,(0,0)}$ & 0.050 & 0.092 & 0.818621 & 0.062 & 0.109 & 0.718379 & 0.052 & 0.104 & 0.612128\\
  &$T_{n,(0,0)}$ & 0.045 & 0.084 & 0.665399 & 0.053 & 0.111 & 0.459543 & 0.045 & 0.104 & 0.226206\\
  \cline{2-11}
  &$R_{n,(1,0)}$ & 0.049 & 0.090 & 0.995881 & 0.060 & 0.106 & 0.902243 & 0.051 & 0.116 & 0.612128\\
  &$S_{n,(1,0)}$ & 0.049 & 0.088 & 0.818621 & 0.062 & 0.101 & 0.818621 & 0.051 & 0.106 & 0.459543\\
  &$T_{n,(1,0)}$ & 0.044 & 0.084 & 0.863178 & 0.056 & 0.101 & 0.863178 & 0.052 & 0.104 & 0.459543\\
  \cline{2-11}
&$W_{n}$       & 0.029 & 0.076 & 0.024117 & 0.036 & 0.085 & 0.111356 & 0.038 & 0.088 & 0.129364\\
  \cline{2-11}
 & $T$           & 0.025 & 0.049 & $<$ 2.2e-16 & 0.034 & 0.065 & 1.00e-07 & 0.024 & 0.058 & 5.30e-06 \\
  &$I_B$         & 0.116 & 0.140 & $<$ 2.2e-16 & 0.141 & 0.162 & $<$ 2.2e-16 & 0.129 & 0.153 & $<$ 2.2e-16 \\
  &$NI_B$        & 0.045 & 0.074 & 6.10e-06 & 0.033 & 0.081 & $<$ 2.2e-16 & 0.029 & 0.063 & $<$ 2.2e-16 \\
  \hline
\end{tabular}} \end{center}\end{table}

To measure the performance of the considered approximations, we  calculated the $p$-value of the Kolmogorov-Smirnov test statistic of uniformity (KS) for each set of 1,000 values obtained for each test statistic. These values were rounded to 2 decimal places.

Looking at these tables we conclude that the asymptotic approximation to the $p$-values works better for $T$ than for $I_B$ and $NI_B$. Nevertheless, none of them give satisfactory results even for $n=70$. By contrast, the bootstrap provides an accurate approximation of the null distribution of $T_{n,a}$ in all tried cases. As for the choice of $a_1$ and $a_2$, we observe that there is no gain in performance when $a_1 \neq a_2$.

\begin{table}\begin{center}
\caption{Simulation results for the probability of type I error, $\theta_1=1.5, \theta_2=1$.}\label{type I error_2}\vspace{5pt}
{
	\scriptsize
\begin{tabular}{|c|l|cc|c|cc|c|cc|c|}
\cline{3-11}\multicolumn{2}{c|}{}
  &\multicolumn{3}{c|}{$n=30$}& \multicolumn{3}{c|}{$n=50$} & \multicolumn{3}{c|}{$n=70$}\\ \hline   
  $\theta_3;\rho$  & Test
  & f05 & f10 & $KS$ & f05 & f10 & $KS$ & f05 & f10 & $KS$\\
  \hline
0.31; 0.25 &$R_{n,(0,0)}$ & 0.054 & 0.104 & 0.226206 & 0.062 & 0.112 & 0.559560 & 0.052 & 0.111 & 0.069329\\
& $S_{n,(0,0)}$ & 0.054 & 0.103 & 0.718379 & 0.050 & 0.104 & 0.291736 & 0.050 & 0.095 & 0.863178\\
  &$T_{n,(0,0)}$ & 0.056 & 0.091 & 0.718379 & 0.060 & 0.098 & 0.989545 & 0.057 & 0.107 & 0.902243\\
  \cline{2-11}
  &$R_{n,(1,0)}$ & 0.050 & 0.094 & 0.149677 & 0.058 & 0.108 & 0.960002 & 0.051 & 0.107 & 0.013476\\
  &$S_{n,(1,0)}$ & 0.047 & 0.095 & 0.226206 & 0.057 & 0.111 & 0.413150 & 0.053 & 0.097 & 0.508494\\
  &$T_{n,(1,0)}$ & 0.043 & 0.091 & 0.369615 & 0.053 & 0.109 & 0.559560 & 0.058 & 0.115 & 0.559560\\
 \cline{2-11}
  &$R_{n,(0,1)}$ & 0.050 & 0.106 & 0.459543 & 0.060 & 0.112 & 0.902243 & 0.056 & 0.109 & 0.016427\\
  &$S_{n,(0,1)}$ & 0.051 & 0.105 & 0.459543 & 0.057 & 0.105 & 0.413150 & 0.056 & 0.099 & 0.508494\\
  &$T_{n,(0,1)}$ & 0.052 & 0.096 & 0.863178 & 0.060 & 0.098 & 0.863178 & 0.056 & 0.106 & 0.718379\\
\cline{2-11}
&$W_{n}$       & 0.022 & 0.066 & 0.041633 & 0.036 & 0.076 & 0.111356 & 0.037 & 0.082 & 0.111356\\
 \cline{2-11}
 & $T$           & 0.018 & 0.046 & 1.00e-07 & 0.021 & 0.060 & 0.000318 & 0.022 & 0.064 & 0.009785 \\
  &$I_B$         & 0.031 & 0.060 & $<$ 2.2e-16 & 0.013 & 0.028 & $<$ 2.2e-16 & 0.007 & 0.014 & $<$ 2.2e-16 \\
 & $NI_B$        & 0.016 & 0.041 & $<$ 2.2e-16 & 0.010 & 0.018 & $<$ 2.2e-16 & 0.004 & 0.009 & $<$ 2.2e-16 \\
  \hline  \hline
 0.62; 0.51 &$R_{n,(0,0)}$ & 0.047 & 0.095 & 0.459543 & 0.045 & 0.095 & 0.863178 & 0.052 & 0.114 & 0.718379\\
 & $S_{n,(0,0)}$ & 0.048 & 0.104 & 0.818621 & 0.049 & 0.091 & 0.818621 & 0.048 & 0.093 & 0.459543\\
 & $T_{n,(0,0)}$ & 0.042 & 0.093 & 0.934732 & 0.043 & 0.098 & 0.934732 & 0.045 & 0.099 & 0.612128\\
  \cline{2-11}
&  $R_{n,(1,0)}$ & 0.045 & 0.095 & 0.863178 & 0.044 & 0.096 & 0.978036 & 0.056 & 0.101 & 0.559560\\
&  $S_{n,(1,0)}$ & 0.051 & 0.088 & 0.459543 & 0.045 & 0.086 & 0.718379 & 0.048 & 0.100 & 0.459543\\
&  $T_{n,(1,0)}$ & 0.037 & 0.096 & 0.769894 & 0.041 & 0.088 & 0.291736 & 0.047 & 0.109 & 0.226206\\
  \cline{2-11}
&  $R_{n,(0,1)}$ & 0.049 & 0.097 & 0.413150 & 0.045 & 0.101 & 0.718379 & 0.054 & 0.104 & 0.902243\\
&  $S_{n,(0,1)}$ & 0.052 & 0.098 & 0.508494 & 0.042 & 0.098 & 0.612128 & 0.051 & 0.091 & 0.413150\\
&  $T_{n,(0,1)}$ & 0.043 & 0.089 & 0.902243 & 0.049 & 0.095 & 0.508494 & 0.046 & 0.088 & 0.329116\\
\cline{2-11}
&$W_{n}$       & 0.026 & 0.055 & 0.003013 & 0.037 & 0.071 & 0.111356 & 0.039 & 0.079 & 0.111356\\
\cline{2-11}
&  $T$           & 0.056 & 0.088 & 0.000526 & 0.050 & 0.104 & 0.011917 & 0.049 & 0.096 & 0.001109 \\
&  $I_B$         & 0.147 & 0.201 & $<$ 2.2e-16 & 0.169 & 0.223 & $<$ 2.2e-16 & 0.147 & 0.196 & $<$ 2.2e-16 \\
&  $NI_B$        & 0.094 & 0.152 & 0.000622 & 0.082 & 0.145 & 0.006666 & 0.076 & 0.120 & 0.078967 \\
  \hline  \hline
0.92; 0.75&  $R_{n,(0,0)}$ & 0.057 & 0.102 & 0.612128 & 0.054 & 0.097 & 0.413150 & 0.046 & 0.090 & 0.863178\\
&  $S_{n,(0,0)}$ & 0.052 & 0.108 & 0.413150 & 0.050 & 0.091 & 0.769894 & 0.044 & 0.094 & 0.559560\\
 & $T_{n,(0,0)}$ & 0.043 & 0.098 & 0.934732 & 0.056 & 0.102 & 0.226206 & 0.041 & 0.085 & 0.413150\\
  \cline{2-11}
  &$R_{n,(1,0)}$ & 0.053 & 0.104 & 0.559560 & 0.055 & 0.103 & 0.508494 & 0.043 & 0.088 & 0.863178\\
  &$S_{n,(1,0)}$ & 0.049 & 0.103 & 0.769894 & 0.050 & 0.093 & 0.612128 & 0.045 & 0.091 & 0.612128\\
  &$T_{n,(1,0)}$ & 0.040 & 0.107 & 0.665399 & 0.052 & 0.110 & 0.172476 & 0.037 & 0.084 & 0.459543\\
  \cline{2-11}
  &$R_{n,(0,1)}$ & 0.055 & 0.110 & 0.459543 & 0.050 & 0.094 & 0.257432 & 0.044 & 0.082 & 0.818621\\
  &$S_{n,(0,1)}$ & 0.051 & 0.108 & 0.665399 & 0.048 & 0.087 & 0.369615 & 0.045 & 0.091 & 0.508494\\
  &$T_{n,(0,1)}$ & 0.046 & 0.089 & 0.769894 & 0.055 & 0.109 & 0.508494 & 0.043 & 0.088 & 0.665399\\
\cline{2-11}
&$W_{n}$       & 0.037 & 0.081 & 0.000714 & 0.042 & 0.079 & 0.111356 & 0.037 & 0.083 & 0.149677\\
  \cline{2-11}
 & $T$           & 0.029 & 0.059 & 1.70e-06 & 0.057 & 0.094 & 0.008821 & 0.078 & 0.109 & 0.065401\\
  &$I_B$         & 0.091 & 0.116 & $<$ 2.2e-16 & 0.209 & 0.239 & $<$ 2.2e-16 & 0.196 & 0.220 & $<$ 2.2e-16\\
 & $NI_B$        & 0.021 & 0.051 & $<$ 2.2e-16 & 0.089 & 0.152 & 0.001554 & 0.094 & 0.149 & 0.003483\\
  \hline
\end{tabular}} \end{center}\end{table}

To study the power we repeated the above experiment for samples with size $n=50$ and we use the same alternative distributions used in NJ (2014), some of which have also been taken as alternatives by other researchers (see, e.g. Loukas and Kemp, 1986; Rayner and Best, 1995, and NJ, 2016).

The parameters of these alternatives were chosen for the same reason given by NJ (2014). We took $a_1=a_2=0$ because, as observed from the results in the previous experiment, there is no gain in performance when $a_1 \neq a_2$ when approximating the probability of type I error. In addition, taking $a_1=a_2=0$ is less time consuming.

 Table \ref{power} displays the alternatives considered and the estimated power for nominal significance level $\alpha=0.05$. The results presented in this table allow us to conclude that the new test proposed in this paper is able to detect all the alternatives treated and with a power as good or better than the other tests based on the bootstrap method, while the non-consistent tests are not able to detect most of these alternatives, especially tests $I_B$ and $NI_B$.

\begin{landscape}
\begin{table}
\caption{Simulation results for the power ($n=50$).}\label{power}\vspace{5pt}
\begin{center}{\footnotesize
\begin{tabular}{|l|ccc|ccccccc|}
\hline
Alternative & $\frac{var(X_1)}{E(X_1)}$ & $\frac{var(X_2)}{E(X_2)}$& $\rho$ & $R_{n,(0,0)}$ & $S_{n,(0,0)}$ & $W_n$ & $T_{n,(0,0)}$ & $T$ & $I_B$ &  $NI_B$\\
\hline\hline
  $BB(1;0.41,0.02,0.01)$ & 0.590 & 0.980 & 0.026 & 0.860 & 0.871 & 0.829 & 0.857 & 0.103 & 0.000 & 0.000\\
  $BB(1;0.41,0.03,0.02)$ & 0.590 & 0.970 & 0.092 & 0.859 & 0.879 & 0.779 & 0.893 & 0.122 & 0.000 & 0.000\\
  $BB(2;0.42,0.02,0.01)$ & 0.580 & 0.980 & 0.023 & 0.726 & 0.677 & 0.682 & 0.746 & 0.251 & 0.005 & 0.005\\
  $BB(2;0.51,0.01,0.01)$ & 0.490 & 0.990 & 0.099 & 0.900 & 0.862 & 0.847 & 0.887 & 0.656 & 0.001 & 0.001\\
  $BB(2;0.61,0.01,0.01)$ & 0.390 & 0.990 & 0.080 & 0.974 & 0.946 & 0.948 & 0.987 & 0.938 & 0.000 & 0.000\\
  \hline\hline
  $BNB(4;0.93,0.01,0.01)$  & 1.930 & 1.010 & 0.143 & 0.793 & 0.793 & 0.809 & 0.851 & 0.853 & 0.860 & 0.853\\
  $BNB(4;0.97,0.01,0.01)$  & 1.970 & 1.010 & 0.141 & 0.815 & 0.815 & 0.802 & 0.912 & 0.872 & 0.880 & 0.864\\
  $BNB(2;0.97,0.97,0.01)$  & 1.970 & 1.970 & 0.493 & 0.938 & 0.908 & 0.891 & 0.941 & 0.895 & 0.629 & 0.987\\
  $BNB(4;0.98,0.01,0.01)$  & 1.980 & 1.010 & 0.141 & 0.832 & 0.830 & 0.846 & 0.925 & 0.876 & 0.889 & 0.873\\
  $BNB(4,0.99,0.01,0.01)$  & 1.990 & 1.010 & 0.140 & 0.823 & 0.812 & 0.817 & 0.880 & 0.878 & 0.881 & 0.878\\
\hline\hline
  $BPP(0.40;(0.2,0.2,0.1);(1.0,0.9,0.1))$ & 1.226 & 1.190 & 0.413 & 0.956 & 0.930 & 0.950 & 0.989 & 0.803 & 0.000 & 0.000\\
  $BPP(0.40;(0.2,0.3,0.1);(0.9,0.9,0.1))$ & 1.190 & 1.131 & 0.361 & 0.932 & 0.895 & 0.913 & 0.927 & 0.747 & 0.000 & 0.000 \\
  $BPP(0.40;(0.8,0.8,0.1);(0.9,1.0,0.4))$ & 1.003 & 1.010 & 0.322 & 0.867 & 0.821 & 0.834 & 0.864 & 0.617 & 0.000 & 0.000 \\
  $BPP(0.45;(0.8,0.8,0.1);(0.9,0.9,0.2))$ & 1.003 & 1.003 & 0.186 & 0.873 & 0.811 & 0.821 & 0.898 & 0.614 & 0.000 & 0.000 \\
  $BPP(0.7;(0.8,0.8,0.1);(0.9,1.1,0.3))$ & 1.003 & 1.021 & 0.208 & 0.864 & 0.809 & 0.893 & 0.941 & 0.600 & 0.000 & 0.000 \\
  \hline\hline
  $BNTA(0.15;0.01,0.01,0.97)$ & 1.990 & 1.990 & 0.995 & 0.800 & 0.802 & 0.835 & 0.898 & 0.615 & 0.003 & 0.682\\
  $BNTA(0.42;0.01,0.01,0.98)$ & 1.990 & 1.990 & 0.995 & 0.908 & 0.896 & 0.907 & 0.949 & 0.665 & 0.003 & 0.849\\
  $BNTA(0.50;0.01,0.01,0.98)$ & 1.990 & 1.990 & 0.995 & 0.925 & 0.919 & 0.831 & 0.921 & 0.684 & 0.000 & 0.888\\
  $BNTA(0.70;0.01,0.01,0.98)$ & 1.990 & 1.990 & 0.995 & 0.937 & 0.919 & 0.777 & 0.923 & 0.730 & 0.001 & 0.899\\
  $BNTA(0.75;0.01,0.01,0.98)$ & 1.990 & 1.990 & 0.995 & 0.932 & 0.922 & 0.796 & 0.978 & 0.717 & 0.001 & 0.910\\
  \hline\hline
  $BLS(0.01,0.01,0.07)$  & 0.156 & 0.156 & 0.197 & 0.876 & 0.930 & 0.902 & 1.000 & 0.800 & 0.000 & 0.000 \\
  $BLS(0.01,0.01,0.25)$  & 0.224 & 0.224 & 0.829 & 0.809 & 0.895 & 0.916 & 0.982 & 0.749 & 0.015 & 0.086 \\
  $BLS(0.26,0.01,0.04)$  & 0.263 & 0.877 & 0.054 & 0.690 & 0.779 & 0.863 & 1.000 & 0.868 & 0.001 & 0.001 \\
  $BLS(3d/7,2d/7,2d/7)^*$ & 1.000 & 1.000 & 0.447 & 0.762 & 0.876 & 0.872 & 0.930 & 0.198 & 0.159 & 0.144\\
  $BLS(3d/4,d/8,d/8)^*$   & 1.000 & 1.000 & 0.267 & 0.942 & 1.000 & 0.981 & 0.909 & 0.249 & 0.205 & 0.191\\
  \hline
  \multicolumn{8}{l}{\normalsize $^*\, d=1-\exp(-1)\approx 0.63212$.}
\end{tabular}}
\end{center}
\end{table}
\end{landscape}

As we stated, the test we propose is faster than its competitors, the Table \ref{CPU} presents the results obtained.

\begin{table}[h]
\caption{Average CPU time (in seconds).}\label{CPU}
\begin{center}
\begin{tabular}{|l|r|r|r|}
\cline{2-4}\multicolumn{1}{c|}{} & \multicolumn{1}{c|}{$n=30$} & \multicolumn{1}{c|}{$n=50$} & \multicolumn{1}{c|}{$n=70$} \\
  \hline
  $R_{n,(0,0)}$ & 40,804.73 & 43,974.45 & 49,328.93 \\
  $S_{n,(0,0)}$ &  3,040.57 &  7,375.56 & 14,502.74 \\
  $W_{n}$       &  1,452.07 &  1,807.28 &  2,142.86 \\  
  $T_{n,(0,0)}$ &    252.31 &    518.03 &    723.42 \\  \hline
\end{tabular}
\end{center}
\end{table}

\subsection{Real data sets}

To end this section, $T_{n,a}$ is applied to a real data set. This data set was analyzed in Berm\'udez (2009), who used two variables, the number of claims for third-party liability ($X_1$) and the number of claims for the rest of guarantees ($X_2$). The original sample comprised a ten percent sample of the automobile portfolio of a major insurance company operating in Spain in 1995. The author assumed that $(X_1,X_2)$ has a BPD, but according to the report shown in Table \ref{Resultados-datos-reales}, the data set is not well modeled by a BPD. The blanks are due to the fact that $W_n$ does not depend on the value of $(a_1, a_2)$.

\begin{table}
	\caption{Results for the real data set $(n = 80,994)$.}\label{Resultados-datos-reales}
	\begin{center}	
		\begin{tabular}{|c|cccc|}
			\cline{2-5} \multicolumn{1}{c|}{}
			&   \multicolumn{4}{c|}{Claims} \\ \hline
			$(a_1,a_2)$ & $R_{n,(a_1,a_2)}$ & $S_{n,(a_1,a_2)}$ & $T_{n,(a_1,a_2)}$ & $W_n$  \\
			\hline
			$(0,0)$  & 0.001 & 0.001 & 0.000 & 0.020\\
			$(1,0)$  & 0.005 & 0.005 & 0.008 &  \\
			$(0,1)$  & 0.003 & 0.006 & 0.010 & \\
			\hline
			\hline
			\multicolumn{1}{|c|}{$\widehat{\theta}_n$} & \multicolumn{4}{c|}{(0.06702119,\,0.08841783,\,0.01394778)}  \\
			\hline
		\end{tabular}
	\end{center}
\end{table}

\subsection{Case $\boldsymbol{\theta_3=0}$} \label{theta3iguala0}
This case has been excluded from $H_0$ because it is a boundary point. 
This situation occurs when the variables $X_1$ and $X_2$ are independent and is analyzed in NJ (2016) where different ways of approaching it are given, besides references are cited for a detailed treatment, even it is a subject for a future research.

\section{Extension of $\boldsymbol{T_{n,w}}$} \label{Extension T}
In principle, the approach can be generalized to the case $m\geq 3$ and we refer to a manuscript that is uploaded to arXive math.\\

To illustrate this situation we will present the case for $m = 3$, in which we need to satisfy $7 = \sum_{i=1}^{3}\binom{3}{i}$ equations to obtain a characterization of the respective Poisson distribution. It can be seen that the number of equations grows following a sum of combinatorial numbers due to differential equations of different order that must be verified, which range from order 1 to order $m$. These equations arise due to the philosophy of the method to characterize the respective Poisson distribution.

\subsection{Trivariate case} \label{Trivariate case}
For this particular case, from section \ref{The general m-variate case}, for $m = 3$, let
\[X_1=Y_1+Y_4, \,\, X_2=Y_2+Y_4,\,\, X_3=Y_3+Y_4,\]
where $Y_1, Y_2, Y_3, Y_4$ are independent Poisson random variables
with means $\theta'_1=\theta_1-\theta_4>0, \theta'_2=\theta_2-\theta_4>0, \theta'_3=\theta_3-\theta_4>0$ and  $\theta_4>0$, respectively. The joint distribution of the vector $(X_1, X_2, X_3)$ is called  a trivariate Poisson distribution (TPD) with parameter $\theta=(\theta_1,\theta_2,\theta_3,\theta_4)$ (see, e.g. Johnson, Kotz and Balakrishnan, 1997; Loukas and Papageorgiou, 1991). The joint pgf of
$(X_1, X_2, X_3)$ is
\begin{equation}\label{fgp-trivariate}
g(u;\theta)=\exp\bigl\{\theta_1(u_1-1)+ \theta_2(u_2-1)+ \theta_3(u_3-1)+\theta_4(u_1 u_2 u_3-u_1-u_2-u_3+2)\bigr\}.
\end{equation}
The empirical counterpart of pgf is epgf of the data given by 
\[g_n(u)=\frac{1}{n}\sum_{i=1}^n u_1^{X_{i1}}u_2^{X_{i2}}u_3^{X_{i3}},\ \, u=(u_1,u_2,u_3).\]
Now, the objective is to test the hypothesis
\[H_{03}:(X_1,  X_2, X_3)\ \mbox{has a trivariate Poisson distribution}.\]

To achieve this new objective, we give a characterization for the TPD.
\begin{proposition}\label{Soluc-EDP-marginales-conjuntas} Let $g(u_1,u_2,u_3;\theta)\,$ be as defined in (\ref{fgp-trivariate}). Then $g(u_1,u_2,u_3;\theta)\,$ is the only pgf satisfying the following system of partial differential equations
\[	\begin{array}{rcl}
	D_1(u;\theta) &=&\frac{\partial g(u_1,1,1)}{\partial u_1}-\theta_1\, g(u_1,1,1)=0,\\[.1 cm]
	D_2(u;\theta)&=&\frac{\partial g(1,u_2,1)}{\partial u_2}-\theta_2\, g(1,u_2,1)=0,\\[.1 cm]
	D_3(u;\theta) &=&\frac{\partial g(1,1,u_3)}{\partial u_3}-\theta_3\, g(1,1,u_3)=0,\\[.1 cm]
	D_4(u;\theta)&=&\frac{\partial^2 g(u)}{\partial u_1 \,\partial u_2}-g(u)[ \{\theta_1+\theta_4(u_2u_3-1)\}\{\theta_2+\theta_4(u_1u_3-1)\}+\theta_4u_3] =0,\\[.1 cm] 
	D_5(u;\theta)&=&\frac{\partial^2 g(u)}{\partial u_1 \,\partial u_3}-g(u)[ \{\theta_1+\theta_4(u_2u_3-1)\}\{\theta_3+\theta_4(u_1u_2-1)\}+\theta_4u_2] =0,\\[.1 cm] 
	D_6(u;\theta)&=&\frac{\partial^2 g(u)}{\partial u_2 \,\partial u_3}-g(u)[ \{\theta_2+\theta_4(u_1u_3-1)\}\{\theta_3+\theta_4(u_1u_2-1)\}+\theta_4u_1] =0,\\[.1 cm] 
	D_7(u;\theta)&=&\frac{\partial^3 g(u)}{\partial u_1 \,\partial u_2 \,\partial u_3}-g(u)h(u_1,u_2,u_3;\theta) =0, \end{array}
\]
	%
	where $\,h(u;\theta)=\prod_{i=1}^{3} \left\{\theta_i+\theta_4\left(\prod_{j\neq i}u_j-1\right)\right\}+\theta_4\left(1+\sum_{k=1}^{3}u_k\left\{\theta_k+\theta_4\left(\prod_{j\neq k}u_j-1\right)\right\}\right)$.
\end{proposition}
By Proposition 1 in NJ (2014), $g(u)$ and its derivatives can be consistently estimated by the epgf and the derivatives of the epgf, respectively. Thus, if $H_{03}$ is true, then the functions
\[
\begin{array}{rcl}
D_{1n}(u;\widehat{\theta}_n) &=&\frac{\partial g_n(u_1,1,1)}{\partial u_1}-\widehat{\theta}_{1n}\, g_n(u_1,1,1)=0,\\[.1 cm]
D_{2n}(u;\widehat{\theta}_n)&=&\frac{\partial g_n(1,u_2,1)}{\partial u_2}-\widehat{\theta}_{2n}\, g_n(1,u_2,1)=0,\\[.1 cm]
D_{3n}(u;\widehat{\theta}_n) &=&\frac{\partial g_n(1,1,u_3)}{\partial u_3}-\widehat{\theta}_{3n}\, g_n(1,1,u_3)=0,\\[.1 cm]
D_{4n}(u;\widehat{\theta}_n)&=&\frac{\partial^2 g_n(u)}{\partial u_1 \,\partial u_2}-g_n(u)[ \{\widehat{\theta}_{1n}+\widehat{\theta}_{4n}(u_2u_3-1)\}\{\widehat{\theta}_{2n}+\widehat{\theta}_{4n}(u_1u_3-1)\}+\widehat{\theta}_{4n}u_3] =0,\\[.1 cm] 
D_{5n}(u;\widehat{\theta}_n)&=&\frac{\partial^2 g_n(u)}{\partial u_1 \,\partial u_3}-g_n(u)[ \{\widehat{\theta}_{1n}+\widehat{\theta}_{4n}(u_2u_3-1)\}\{\widehat{\theta}_{3n}+\widehat{\theta}_{4n}(u_1u_2-1)\}+\widehat{\theta}_{4n}u_2] =0,\\[.1 cm] 
D_{6n}(u;\widehat{\theta}_n)&=&\frac{\partial^2 g_n(u)}{\partial u_2 \,\partial u_3}-g_n(u)[ \{\widehat{\theta}_{2n}+\widehat{\theta}_{4n}(u_1u_3-1)\}\{\widehat{\theta}_{3n}+\widehat{\theta}_{4n}(u_1u_2-1)\}+\widehat{\theta}_{4n}u_1] =0,\\[.1 cm] 
D_{7n}(u;\widehat{\theta}_n)&=&\frac{\partial^3 g_n(u)}{\partial u_1 \,\partial u_2 \,\partial u_3}-g_n(u)h(u_1,u_2,u_3;\widehat{\theta}_n) =0, \end{array}
\]
should be close to $0$, $\forall (u_1,u_2,u_3) \in [0,1]^3$, where $\widehat{\theta}_n=(\widehat{\theta}_{1n}, \widehat{\theta}_{2n}, \widehat{\theta}_{3n})$ is a consistent estimator of $\theta$ and
$g_n(u_1,u_2,u_3)$ is the epgf associated with the data, i. e.,
\[g_n(u_1,u_2,u_3)=\frac{1}{n}\sum_{i=1}^n u_1^{X_{i1}}u_2^{X_{i2}}u_3^{X_{i3}}.\]

Thus, to test $H_{03}$ we consider the following test statistic
\begin{equation}\notag
T_{3,n,w}=n\int_0^1\int_0^1\int_0^1 \left\{D^2_{1n}\bigl(u;\widehat{\theta}_n \bigr)+D^2_{2n}\bigl(u;\widehat{\theta}_{n} \bigr)+\cdots+D^2_{7n}\bigl(u;\widehat{\theta}_{n} \bigr)\right\}w(u)\ du,
\end{equation}
where $w(u)$ is a measurable non-negative function with finite integral over $[0,1]^3$. Similar results to those stated in  Sections \ref{A characterization of the BPD}, \ref{Null bootstrap distribution}, and \ref{Alternatives} for the bivariate case can be established for $T_{3,n,w}$.

\begin{remark}
	So far we have not managed to obtain numerical results for the case $m = 3$ due to the large number of calculations involved in $T_{3, n, w}$. We can assure that this new test is not recommended for $m \geq 3$ and it is preferable to use the $W_n$ statistic.
\end{remark}

\subsection{Simulated data for the trivariate case}
To simulate type I error we follow a procedure similar to that described for the bivariate case, but we do not have competitors. We consider three situations: a) $\theta_1 = \theta_2 = \theta_3$, b) $\theta_1 = \theta_2 \neq \theta_3$, $\theta_2 = \theta_3 \neq \theta_1$ and c) $\theta_1 \neq \theta_2$ and $\theta_1\neq \theta_3$ and $\theta_2\neq \theta_3$. In each of these cases $\theta_1$, $\theta_2$, $\theta_3> \theta_4$. In addition, $\theta_4$ was chosen in such a way that the correlation coefficients, $\rho=(\rho_{12}, \rho_{13}, \rho_{23})$, were equal or very close to 0.25, 0.5, 0.75 and 1.00,  where $\rho_{ij}=\frac{Cov(X_i,X_j)}{\sqrt{Var(X_i)Var(X_j)}}$.\\

Tables \ref{type I error_3} and \ref{type I error_4} display the fraction of estimated $p$-values less than or equal to $0.05$ and $0.10$, which are the estimated type I error probabilities for $\alpha=0.05$ and $0.10$ (denoted as f05 and f10 in the tables), respectively.\\

\begin{table}\begin{center}
		\caption{Simulation results for the probability of type I error, $\theta_1=\theta_2=\theta_3=1$.}\label{type I error_3}
		{\footnotesize
			\begin{tabular}{|c|l|cc|c|cc|c|cc|c|}
				\cline{3-11}\multicolumn{2}{c|}{}
				&\multicolumn{3}{c|}{$n=30$}& \multicolumn{3}{c|}{$n=50$} & \multicolumn{3}{c|}{$n=70$}\\ 
				\hline 
				$\theta_4=\rho_0$ $^{(*)}$ & Test & f05 & f10 & $KS$ & f05 & f10 & $KS$ & f05 & f10 & $KS$\\
				\hline
				0.25 &$R_{3,n,(0,0,0)}$ & 0.035 & 0.083 & 0.329116 & 0.044 & 0.089 & 0.129364 & 0.047 & 0.092 & 0.934732\\
				&$S_{3,n,(0,0,0)}$ & 0.040 & 0.087 & 0.508494 & 0.043 & 0.091 & 0.818621 & 0.048 & 0.093 & 0.863178\\
				&$T_{3,n,(0,0,0)}$ & 0.041 & 0.088 & 0.111356 & 0.043 & 0.090 & 0.508494 & 0.048 & 0.096 & 0.718379\\
				\cline{2-11}
				&$R_{3,n,(1,0,0)}$ & 0.039 & 0.089 & 0.329116 & 0.042 & 0.091 & 0.329116 & 0.045 & 0.094 & 0.818621\\
				&$S_{3,n,(1,0,0)}$ & 0.038 & 0.090 & 0.902243 & 0.043 & 0.090 & 0.508494 & 0.043 & 0.090 & 0.960002\\
				&$T_{3,n,(1,0,0)}$ & 0.040 & 0.088 & 0.257432 & 0.041 & 0.090 & 0.329116 & 0.046 & 0.094 & 0.369615\\
				\cline{2-11}
				&$W_{3,n}$       & 0.032 & 0.075 & 0.129364 & 0.038 & 0.082 & 0.129364 & 0.042 & 0.090 & 0.612128\\
				\hline
				\hline
				0.75 & $R_{3,n,(0,0,0)}$ & 0.038 & 0.82 & 0.049545 & 0.041 & 0.091 & 0.197933 & 0.042 & 0.093 & 0.559560\\
				&$S_{3,n,(0,0,0)}$ & 0.040 & 0.090 & 0.069329 & 0.042 & 0.092 & 0.257432 & 0.047 & 0.094 & 0.413150\\
				&$T_{3,n,(0,0,0)}$ & 0.040 & 0.088 & 0.149677 & 0.042 & 0.085 & 0.257432 & 0.051 & 0.101 & 0.413150\\
				\cline{2-11}
				&$R_{3,n,(1,0,0)}$ & 0.041 & 0.081 & 0.129364 & 0.044 & 0.087 & 0.172476 & 0.046 & 0.102 & 0.863178\\
				&$S_{3,n,(1,0,0)}$ & 0.040 & 0.088 & 0.129364 & 0.045 & 0.090 & 0.197933 & 0.046 & 0.095 & 0.291736\\
				&$T_{3,n,(1,0,0)}$ & 0.041 & 0.090 & 0.111356 & 0.043 & 0.094 & 0.413150 & 0.051 & 0.112 & 0.459543\\
				\cline{2-11}
				&$W_{3,n}$       & 0.032 & 0.081 & 0.013476 & 0.042 & 0.087 & 0.111356 & 0.043 & 0.091 & 0.129364\\
				\hline
				\multicolumn{3}{l}{$^{(*)} \rho_{12}=\rho_{13}=\rho_{23}=\rho_0$}
\end{tabular}} \end{center}\end{table}

\begin{table}\begin{center}
		\caption{Simulation results for the probability of type I error, $\theta_1=\theta_2=\theta_3=2$.}\label{type I error_4}
		{\footnotesize
			\begin{tabular}{|c|l|cc|c|cc|c|cc|c|}
				\cline{3-11}\multicolumn{2}{c|}{}
				&\multicolumn{3}{c|}{$n=30$}& \multicolumn{3}{c|}{$n=50$} & \multicolumn{3}{c|}{$n=70$}\\ 
				\hline 
				$\theta_4=2\rho_0$ $^{(*)}$& Test & f05 & f10 & $KS$ & f05 & f10 & $KS$ & f05 & f10 & $KS$\\
				\hline
				0.5 &$R_{3,n,(0,0,0)}$ & 0.039 & 0.085 & 0.111356 & 0.042 & 0.091 & 0.257432 & 0.043 & 0.092 & 0.863178\\
				&$S_{3,n,(0,0,0)}$ & 0.040 & 0.082 & 0.508494 & 0.042 & 0.088 & 0.818621 & 0.046 & 0.093 & 0.934732\\
				&$T_{3,n,(0,0,0)}$ & 0.041 & 0.085 & 0.329116 & 0.042 & 0.091 & 0.508494 & 0.045 & 0.094 & 0.718379\\
				\cline{2-11}
				&$R_{3,n,(1,0,0)}$ & 0.038 & 0.087 & 0.329116 & 0.041 & 0.093 & 0.329116 & 0.048 & 0.094 & 0.818621\\
				&$S_{3,n,(1,0,0)}$ & 0.039 & 0.089 & 0.902243 & 0.043 & 0.092 & 0.508494 & 0.043 & 0.092 & 0.960002\\
				&$T_{3,n,(1,0,0)}$ & 0.040 & 0.090 & 0.129364 & 0.042 & 0.090 & 0.329116 & 0.045 & 0.095 & 0.369615\\
				\cline{2-11}
				&$W_{3,n}$       & 0.032 & 0.086 & 0.111356 & 0.043 & 0.087 & 0.369615 & 0.048 & 0.092 & 0.612128\\
				\hline
				\hline
				1.00 & $R_{3,n,(0,0,0)}$ & 0.038 & 0.091 & 0.129364 & 0.044 & 0.092 & 0.197933 & 0.046 & 0.094 & 0.559560\\
				&$S_{3,n,(0,0,0)}$ & 0.042 & 0.089 & 0.049545 & 0.043 & 0.093 & 0.257432 & 0.046 & 0.095 & 0.413150\\
				&$T_{3,n,(0,0,0)}$ & 0.040 & 0.090 & 0.149677 & 0.043 & 0.091 & 0.257432 & 0.053 & 0.101 & 0.413150\\
				\cline{2-11}
				&$R_{3,n,(1,0,0)}$ & 0.051 & 0.101 & 0.129364 & 0.044 & 0.097 & 0.172476 & 0.047 & 0.109 & 0.863178\\
				&$S_{3,n,(1,0,0)}$ & 0.040 & 0.090 & 0.069329 & 0.056 & 0.109 & 0.197933 & 0.052 & 0.102 & 0.291736\\
				&$T_{3,n,(1,0,0)}$ & 0.040 & 0.089 & 0.129364 & 0.044 & 0.094 & 0.413150 & 0.051 & 0.103 & 0.459543\\
				\cline{2-11}
				&$W_{3,n}$       & 0.038 & 0.081 & 0.111356 & 0.042 & 0.087 & 0.111356 & 0.047 & 0.092 & 0.413150\\
				\hline
				\multicolumn{3}{l}{$^{(*)} \rho_{12}=\rho_{13}=\rho_{23}=\rho_0$}
\end{tabular}} \end{center}\end{table}

\begin{table}\begin{center}
		\caption{Simulation results for the probability of type I error, $\theta_1=\theta_2=0.2$.}\label{type I error_5}
		{
			\scriptsize
			\begin{tabular}{|c|l|cc|c|cc|c|cc|c|}
				\cline{3-11}\multicolumn{2}{c|}{}
				&\multicolumn{3}{c|}{$n=30$}& \multicolumn{3}{c|}{$n=50$} & \multicolumn{3}{c|}{$n=70$}\\ \hline   
				$\theta_3, \theta_4; \rho$	& Test
				& f05 & f10 & $KS$ & f05 &  f10 & $KS$ & f05 & f10 & $KS$\\
				\hline
				$0.8, 0.1; (0.5,0.25,0.25)$&$R_{3,n,(0,0,0)}$ & 0.046 & 0.094 & 0.069329 & 0.058 & 0.107 & 0.559560 & 0.052 & 0.101 & 0.226206\\
				 & $S_{3,n,(0,0,0)}$ & 0.045 & 0.087 & 0.291736 & 0.050 & 0.104 & 0.718379 & 0.051 & 0.093 & 0.863178\\
				&$T_{3,n,(0,0,0)}$ & 0.046 & 0.090 & 0.718379 & 0.056 & 0.098 & 0.902243 & 0.053 & 0.105 & 0.989545\\
				\cline{2-11}
				&$R_{3,n,(1,0,0)}$ & 0.048 & 0.091 & 0.149677 & 0.054 & 0.108 & 0.013476 & 0.052 & 0.107 & 0.960002\\
				&$S_{3,n,(1,0,0)}$ & 0.043 & 0.092 & 0.226206 & 0.055 & 0.111 & 0.413150 & 0.053 & 0.097 & 0.508494\\
				&$T_{3,n,(1,0,0)}$ & 0.040 & 0.091 & 0.369615 & 0.053 & 0.109 & 0.559560 & 0.052 & 0.105 & 0.559560\\
				\cline{2-11}
				&$W_{3,n}$       & 0.038 & 0.086 & 0.041633 & 0.043 & 0.087 & 0.111356 & 0.047 & 0.092 & 0.197933\\
				\hline 
\end{tabular}} \end{center}\end{table}

\begin{table}\begin{center}
		\caption{Simulation results for the probability of type I error, $\theta_2=\theta_3=0.8$.}\label{type I error_6}
		{
			\scriptsize
			\begin{tabular}{|c|l|cc|c|cc|c|cc|c|}
				\cline{3-11}\multicolumn{2}{c|}{}
				&\multicolumn{3}{c|}{$n=30$}& \multicolumn{3}{c|}{$n=50$} & \multicolumn{3}{c|}{$n=70$}\\ \hline   
				$\theta_1, \theta_4; \rho$	& Test
				& f05 & f10 & $KS$ & f05 &  f10 & $KS$ & f05 & f10 & $KS$\\
				\hline
				$1.8, 0.6; (0.5,0.5,0.75)$&$R_{3,n,(0,0,0)}$ & 0.040 & 0.089 & 0.069329 & 0.045 & 0.090 & 0.226206 & 0.045 & 0.091 & 0.559560\\
				&$S_{3,n,(0,0,0)}$ & 0.042 & 0.086 & 0.291736 & 0.045 & 0.094 & 0.718379 & 0.047 & 0.093 & 0.902243\\
				&$T_{3,n,(0,0,0)}$ & 0.044 & 0.089 & 0.197933 & 0.046 & 0.092 & 0.863178 & 0.047 & 0.101 & 0.989545\\
				\cline{2-11}
				&$R_{3,n,(1,0,0)}$ & 0.046 & 0.089 & 0.013476 & 0.045 & 0.094 & 0.508494 & 0.046 & 0.102 & 0.960002\\
				&$S_{3,n,(1,0,0)}$ & 0.044 & 0.088 & 0.226206 & 0.045 & 0.091 & 0.413150 & 0.045 & 0.098 & 0.559560\\
				&$T_{3,n,(1,0,0)}$ & 0.041 & 0.089 & 0.369615 & 0.045 & 0.093 & 0.508494 & 0.046 & 0.101 & 0.559560\\
				\cline{2-11}
				&$W_{3,n}$        & 0.039 & 0.087 & 0.149677 & 0.043 & 0.087 & 0.197933 & 0.046 & 0.095 & 0.291736\\
				\hline 
\end{tabular}} \end{center}\end{table}

\begin{table}\begin{center}
		\caption{Simulation results for the probability of type I error, $\theta_1=8.7$, $\theta_2=8.8$, $\theta_3=8.9$.}\label{type I error_7}
		{
			\scriptsize
			\begin{tabular}{|c|l|cc|c|cc|c|cc|c|}
				\cline{3-11}\multicolumn{2}{c|}{}
				&\multicolumn{3}{c|}{$n=30$}& \multicolumn{3}{c|}{$n=50$} & \multicolumn{3}{c|}{$n=70$}\\ \hline   
				$\theta_4; \rho$	& Test
				& f05 & f10 & $KS$ & f05 &  f10 & $KS$ & f05 & f10 & $KS$\\
				\hline
				$2.2; (0.251,0.250,0.249)$&$R_{3,n,(0,0,0)}$ & 0.036 & 0.083 & 0.149677 & 0.043 & 0.089 & 0.226206 & 0.041 & 0.089 & 0.863178\\
				&$S_{3,n,(0,0,0)}$ & 0.038 & 0.085 & 0.291736 & 0.040 & 0.089 & 0.508494 & 0.043 & 0.092 & 0.718379\\
				&$T_{3,n,(0,0,0)}$ & 0.040 & 0.086 &0.226206  & 0.042 & 0.089 & 0.559560 & 0.045 & 0.105 & 0.902243\\
				\cline{2-11}
				&$R_{3,n,(1,0,0)}$ & 0.036 & 0.085 & 0.197933 & 0.041 & 0.089 & 0.413150 & 0.045 & 0.104 & 0.863178\\
				&$S_{3,n,(1,0,0)}$ & 0.037 & 0.086 & 0.369615 & 0.040 & 0.090 & 0.508494 & 0.044 & 0.094 & 0.559560\\
				&$T_{3,n,(1,0,0)}$ & 0.034 & 0.084 & 0.197933 & 0.042 & 0.090 & 0.508494 & 0.045 & 0.105 & 0.559560\\
				\cline{2-11}
				&$W_{3,n}$        & 0.035 & 0.083 & 0.069329 & 0.040 & 0.087 & 0.197933 & 0.046 & 0.093 & 0.291736\\
				\hline 
\end{tabular}} \end{center}\end{table}

\begin{table}\begin{center}
		\caption{Simulation results for the probability of type I error, $\theta_1=9.7$, $\theta_2=9.6$, $\theta_3=9.5$.}\label{type I error_8}
		{
			\scriptsize
			\begin{tabular}{|c|l|cc|c|cc|c|cc|c|}
				\cline{3-11}\multicolumn{2}{c|}{}
				&\multicolumn{3}{c|}{$n=30$}& \multicolumn{3}{c|}{$n=50$} & \multicolumn{3}{c|}{$n=70$}\\ \hline   
				$\theta_4; \rho$	& Test
				& f05 & f10 & $KS$ & f05 &  f10 & $KS$ & f05 & f10 & $KS$\\
				\hline
				$2.4; (0.249,0.250,0.251)$&$R_{3,n,(0,0,0)}$ & 0.035 & 0.081 & 0.413150 & 0.041 & 0.086 & 0.508494 & 0.041 & 0.088 & 0.863178\\
				&$S_{3,n,(0,0,0)}$ & 0.036 & 0.086 & 0.291736 & 0.039 & 0.087 & 0.369615 & 0.044 & 0.091 & 0.718379\\
				&$T_{3,n,(0,0,0)}$ & 0.039 & 0.086 & 0.226206 & 0.040 & 0.088 & 0.559560 & 0.044 & 0.095 & 0.508494\\
				\cline{2-11}
				&$R_{3,n,(1,0,0)}$ & 0.037 & 0.085 & 0.197933 & 0.040 & 0.088 & 0.149677 & 0.046 & 0.094 & 0.559560\\
				&$S_{3,n,(1,0,0)}$ & 0.036 & 0.085 & 0.226206 & 0.041 & 0.089 & 0.291736 & 0.045 & 0.093 & 0.369615\\
				&$T_{3,n,(1,0,0)}$ & 0.035 & 0.084 & 0.069329 & 0.040 & 0.090 & 0.559560 & 0.046 & 0.095 & 0.508494\\
				\cline{2-11}
				&$W_{3,n}$        & 0.034 & 0.082 & 0.197933 & 0.040 & 0.088 & 0.197933 & 0.045 & 0.094 & 0.863178\\
				\hline 
\end{tabular}} \end{center}\end{table}

As we had anticipated in Remark 1, the proposed new  test $T_{3,n}$ is not faster than some of its competitors, as can be seen in the Table \ref{CPU_2}.

\begin{table}
	\caption{Average CPU time (in seconds).}\label{CPU_2}
	\begin{center}
		\begin{tabular}{|l|r|r|r|}
			\cline{2-4}\multicolumn{1}{c|}{} & \multicolumn{1}{c|}{$n=30$} & \multicolumn{1}{c|}{$n=50$} & \multicolumn{1}{c|}{$n=70$} \\
			\hline
			$R_{3,n,(0,0,0)}$ & 47,338.59 & 72,539.85 & 107,952.20 \\
			$S_{3,n,(0,0,0)}$ &  3,234.27 &  9,357.65 &  18,252.43 \\
			$T_{3,n,(0,0,0)}$ &  4,486.30 & 12,617.74 &  24,763.87 \\
			$W_{3,n}$       &  1,833.80 &  2,174.36 &   2,323.33\\  
		  \hline
		\end{tabular}
	\end{center}
\end{table}

\subsection{Real data set for trivariate case}
The data set was analyzed in Catalina Bolanc\'e \& Raluca Vernic (2017), the data come from the Spanish insurance market and consist of a random sample of 162,019 policyholders who had had one or more auto and home policies during the decade 2006-2015. Catalina Bolanc\'e \& Raluca Vernic (2017) used three dependent variables: the number of claims in auto insurance at fault involving only property damage ($X_1$); the number of claims in auto insurance at fault with bodily injury ($X_2$); and, the number of claims in home insurance at fault ($X_3$).\\

Table \ref{Resultados-datos-reales-trivar} shows the p-values obtained by applying the test we propose. It is concluded that the data do not come from a trivariate Poisson distribution, this is in agreement with the researchers who used this data set to model a trivariate Sarmanov distribution.

\begin{table}
	\caption{Results for the real data set ($n=162,019$).}\label{Resultados-datos-reales-trivar}
	\begin{center}	
		\begin{tabular}{|c|cccc|}
			\cline{2-5} \multicolumn{1}{c|}{}
			&   \multicolumn{4}{c|}{Claims}\\ \hline
			$(a_1,a_2,a_3)$ & $R_{3,n,(a_1,a_2,a_3)}$ & $S_{3,n,(a_1,a_2,a_3)}$ & $T_{3,n,(a_1,a_2,a_3)}$ & $W_{3,n}$ \\
			\hline
			$(0,0,0)$  & 0.001 & 0.001 & 0.001 & 0.0001\\
			$(1,0,0)$  & 0.002 & 0.002 & 0.001 &  \\
			\hline
			\hline
			\multicolumn{1}{|c|}{$\widehat{\theta}_n$} & \multicolumn{4}{c|}{(0.249051, 0.03508231, 0.201069, 0.03508218
				)} \\
			\hline
		\end{tabular}
	\end{center}
\end{table}

\section*{Acknowledgements}
The author would like to thank the Departamento de Investigaci\'on de la Universidad del B\'io-B\'io and the Grupo de Investigaci\'on Matem\'atica Aplicada GI 172409/C de la Universidad del B\'io-B\'io, Chile. He also thanks the anonymous reviewers and the editor of this journal for their valuable time and their careful comments and suggestions with which the quality of this paper has been improved.
Thanks to Florencia Osorio for her valuable comments in the revision of the manuscript.

\section*{Appendix}
\subsection*{Proofs}

The proofs of Theorems 1 to 4 are quite similar to those of Theorems 3, 1, 2 and 4 in NJ (2014), respectively. Here we give a sketch of the proof of Proposition \ref{Soluc-EDP-conjunta-y-marginales}. The proof of Proposition 2 follows the steps of the proof of Proposition 1 by occupying the recurrence relationships for the probabilities and their respective partial derivatives given in Loukas and Papageorgiou (1991). A detailed derivation of the results can be obtained from the authors upon request.\\

\noindent {\bf Proof of Proposition \ref{Soluc-EDP-conjunta-y-marginales}} \hspace{2pt} Let $(X_1,X_2)$ be a random vector and let $g(u_1,u_2)=\sum_{i,j\geq 0}P_{ij}\, u_1^i u_2^j$ be its pgf, where $P_{ij}=P(X_1=i,X_2=j)$. Let $f(u_1,u_2;\theta)=c_0+c_1u_1+c_2u_2+c_3u_1u_2,\,$ with $\,c_0=\theta_3+(\theta_1-\theta_3)(\theta_2-\theta_3),\ c_1=\theta_3(\theta_1-\theta_3),$ $c_2=\theta_3(\theta_2-\theta_3)\,$ and $\, c_3=\theta_3^2.\,$ Then
\begin{align}
\frac{\partial^2 g(u_1,u_2)}{\partial u_1\, \partial u_2}&=\sum_{i,j\geq 1}P_{ij}\, i j\, u_1^{i-1} u_2^{j-1}=\sum_{i,j\geq 0}P_{i+1,j+1}\, (i+1)(j+1)\, u_1^i u_2^j,\notag \\[.1 cm]
f(u_1,u_2;\theta)g(u_1,u_2)&=c_0 P_{00}+\sum_{i\geq 1}\!\left(c_0 P_{i0}+c_1 P_{i-1,0}\right)u_1^i+\sum_{j\geq 1}\!\left(c_0 P_{0j}+c_2 P_{0,j-1}\right)u_2^j\notag \\
&\quad +\sum_{i,j\geq 1}\!\left(c_0 P_{ij}+c_1 P_{i-1,j}+c_2 P_{i,j-1}+c_3 P_{i-1,j-1}\right)u_1^iu_2^j.\notag
\end{align}

From the first equation in (\ref{EDP-conjunta-y-marginales}), $D_1(u;\theta)=0$, then by matching coefficients, we obtain
\begin{equation}\label{ec-recurrencia}
\left.
  \begin{array}{rcl}
P_{11} &=& c_0 P_{00}, \\[.15 cm]
(i+1)P_{i+1,1} &=& c_0 P_{i0}+c_1 P_{i-1,0},\ i\in \mathbb{N}, \\[.15 cm]
 (j+1)P_{1,j+1} &=& c_0 P_{0j}+c_2 P_{0,j-1},\ j\in \mathbb{N}, \\[.15 cm]
(i+1)(j+1)P_{i+1,j+1} &=& c_0 P_{ij}+c_1 P_{i-1,j}+c_2 P_{i,j-1}+c_3 P_{i-1,j-1},\ i,j\in \mathbb{N}.
  \end{array}
  \right\}
\end{equation}

With enough algebraic work we can demonstrate that equations (\ref{ec-recurrencia}) satisfy (1) or (2) and (5) in Kawamura (1985). Moreover, the last two equations in (\ref{EDP-conjunta-y-marginales}) satisfy (3) and (4) in Kawamura (1985). Therefore, the result is obtained by applying Theorem 3 in Kawamura (1985).

\newpage
\section*{BIBLIOGRAPHY}

\noindent Baringhaus, L., and Henze, N. (1992). A goodness of fit test for the Poisson distribution based on the empirical generating function. {\it Statistics $\&$ Probability Letters}, {\bf 13}, 269--274.

\vskip 3mm

\noindent Berm\'udez, L. (2009). A priori ratemaking using bivariate Poisson regression models. {\it Insurance: Mathematics and Economics}, {\bf 44}, 135--141.

\vskip 3mm

\noindent Bolanc\'e, C., and Vernic, R. (2017). Multivariate count data generalized linear models: Three approaches based on the Sarmanov distribution. {\it Research Institute of Applied Economics}, {\bf 1}, 1--25.

\vskip 3mm

\noindent Crockett, N. G. (1979). {\it A quick test of fit of a bivariate distribution}. In Interactive Statistics, D. McNeil (ed.), 185--191. Amsterdam: North-Holland.

\vskip 3mm

\noindent Haight, F. A. (1967). {\it Handbook of the Poisson distribution}. New York: John Wiley \& Sons.

\vskip 3mm

\noindent Holgate, P. (1964). Estimation for the Bivariate Poisson Distribution. {\it Biometrika}, {\bf 51}, 241--245.

\vskip 3mm

\noindent Janssen, A. (2000). Global power functions of goodness of fit tests. {\it The Annals of Statistics}, {\bf 28}, 239--253.

\vskip 3mm

\noindent Johnson, N. L., and Kotz, S. (1969). {\it Distributions in Statistics: Discrete Distributions}. Wiley, New York.

\vskip 3mm

\noindent Johnson, N. L., Kotz, S., and Balakrishnan, N. (1997). {\it Discrete Multivariate Distributions}. Wiley, New York.

\vskip 3mm

\noindent Karlis, D., and Tsiamyrtzis, P. (2008). Exact Bayesian modeling for bivariate Poisson data and extensions. {\it Statistics and Computing}, {\bf 18}, 27--40.

\vskip 3mm

\noindent Kawamura, K. (1985). A note on the recurrent relations for the bivariate Poisson distribution. {\it Kodai Math. J.}, {\bf 8}, 70--78.

\vskip 3mm

\noindent Kocherlakota, S., and Kocherlakota, K. (1992). {\it Bivariate Discrete Distributions}. Marcel Dekker, Inc, New York.

\vskip 3mm

\noindent Loukas, S., and Kemp, C. D. (1986). The Index of Dispersion Test for the Bivariate Poisson Distribution. {\it Biometrics}, {\bf 42}, 941--948.

\vskip 3mm

\noindent Loukas, S., and Papageorgiou, H. (1991). On a trivariate Poisson Distribution. {\it Applications of Mathematics}, {\bf 6}, 432--439.

\vskip 3mm

\noindent Nakamura, M., and P\'erez-Abreu, V. (1993). Use of an Empirical Probability Generating Function for Testing a Poisson Model. {\it Canadian Journal of Statistics}, {\bf 21}, 149--156.

\vskip 3mm

\noindent Novoa-Mu\~noz, F., and Jim\'enez-Gamero, M. D. (2014). Testing for the bivariate Poisson distribution. {\it Metrika}, {\bf 77}, 771--793.

\vskip 3mm

\noindent Novoa-Mu\~noz, F., and Jim\'enez-Gamero, M. D. (2016). A goodness-of-fit test for the multivariate Poisson distribution. {\it SORT}, {\bf 40}, 113--138.

\vskip 3mm

\noindent R Core Team. (2017). {\it R: A Language and Environment for Statistical Computing. R Foundation for Statistical Computing, Vienna, Austria}. http://www.R-project.org.

\vskip 3mm

\noindent Rayner, J. C. W., and Best, D. J. (1995). Smooth Tests for the Bivariate Poisson Distribution. {\it Australian $\&$ New Zealand Journal of Statistics}, {\bf 37}, 233--245.

\end{document}